\documentclass[11pt]{amsart}

\usepackage{amssymb}
\usepackage{latexsym}
\usepackage{amsmath}
\usepackage{amsthm}
\usepackage{amsfonts}
\usepackage{color}
\usepackage{epsfig}
\usepackage{pst-grad} 
\usepackage{pst-plot} 

\newcommand{\pboundary}{\partial^l Q}
\newcommand{\pexterior}{Q^{ext}}
\newcommand{\pboundaryT}{\partial^l Q_T}
\newcommand{\pexteriorT}{Q_T^{ext}}
\newcommand{\R}{\mathbb{R}}
\newcommand{\N}{\mathbb{N}}

\newcommand{\I}{\mathcal{I}}
\newcommand{\h}{\mathbb{H}}

\theoremstyle{plain}
\newtheorem{defi}{Definition}[section]
\newtheorem{prop}[defi]{Proposition}
\newtheorem{teo}[defi]{Theorem}

\newtheorem{lema}[defi]{Lemma}

\newtheorem{remark}[defi]{Remark}

\theoremstyle{definition}

\theoremstyle{remark}

\numberwithin{equation}{section}

\begin{document}

\title[]{Existence, Uniqueness and Asymptotic Behavior for Nonlocal Parabolic Problems with Dominating Gradient Terms.}

\author[]{Guy Barles}
\address{
Guy Barles:
Laboratoire de Math\'ematiques et Physique Th\'eorique (UMR CNRS 7350), F\'ed\'eration Denis Poisson (FR CNRS 2964),
Universit\'e Fran\c{c}ois Rabelais Tours, Parc de Grandmont, 37200 Tours, France. {\tt Guy.Barles@lmpt.univ-tours.fr} }

\author[]{Erwin Topp}
\address{
Erwin Topp:
Departamento de Ingenier\'\i a Matem\'atica (UMI 2807 CNRS), Universidad de Chile, Casilla 170, Correo 3, Santiago, Chile.
and Laboratoire de Math\'ematiques et Physique Th\'eorique (CNRS UMR 7350), Universit\'e Fran\c{c}ois Rabelais,
Parc de Grandmont, 37200, Tours, France. {\tt etopp@dim.uchile.cl} }



\date{\today}

\begin{abstract} In this paper we deal with the well-posedness of Dirichlet problems associated to nonlocal Hamilton-Jacobi parabolic equations in a 
bounded, smooth domain $\Omega$, in the case when the classical boundary condition may be lost. We address the problem for both coercive and 
noncoercive Hamiltonians: for coercive Hamiltonians, our results rely more on the regularity properties of the solutions, while noncoercive case 
are related to optimal control problems and the arguments are based on a careful study of the dynamics near the boundary of the domain. 
Comparison principles for bounded sub and supersolutions are obtained in the context of viscosity solutions with generalized 
boundary conditions, and consequently we obtain the existence and uniqueness of solutions in $C(\bar{\Omega} \times [0,+\infty))$ 
by the application of Perron's method. Finally, we prove that the solution of these problems converges to the solutions of the associated 
stationary problem as $t \to +\infty$ under suitable assumptions on the data.
\end{abstract}

\maketitle

\section{Introduction.}

In this paper we are concerned with the existence, uniqueness and asymptotic 
behavior for the solution of the following Cauchy problem set in $Q = \Omega \times (0,+\infty)$ where $\Omega \subset \R^n$ is a bounded domain 
with smooth boundary
\begin{equation}\tag{CP}\label{cauchy}
\left \{ \begin{array}{rll} \partial_t u 
-\I(u(\cdot, t), x) + H(x, t, u, D u) = &0, \ &\mbox{in} \ Q \\
u(x,t) = & \varphi(x,t), \ &\mbox{in} \ \pexterior \\
u(x,0) = & u_0(x), \ &\mbox{in} \ \bar{\Omega}. \end{array} \right .
\end{equation}
where $u: \R^n \times [0,+\infty) \to \R$ stands for 
the unknown function depending on the ``space'' variable $x \in \R^n$ and the ``time''
variable $t \in [0,+\infty)$, $\partial_t u$ is the derivative of $u$ with respect to $t$ and $Du$ is its gradient with respect to $x$. 
We denote by $\pexterior = \Omega^c \times (0, +\infty)$ and the function $\varphi: \bar{Q}^{ext} \to \R$ is  assumed to be continuous and bounded; 
it represents the prescribed value of $u$ in $\pexterior$ (``Dirichlet boundary condition'').

For $\alpha \in (0,2)$ fixed, $\I$ represents an \textsl{integro-differential operator of order less or equal than $\alpha$}, defined 
in the following way: for $x \in \R^n$ and $\phi$ regular enough at $x$ and bounded in $\R^n$, $\I(\phi, x)$ has the general form
\begin{equation}\label{operator}
\I(\phi, x) = \int_{\R^n} [\phi(x + z) - \phi(x) - \mathbf{1}_B \langle D\phi(x), z \rangle] K(z)|z|^{-(n + \alpha)}dz,
\end{equation}
where $K: \R^n \to \R$ is a measurable, nonnegative and bounded function. Such an operator is called \textsl{elliptic}, and ranges 
from \textsl{zero-th order non local operators} in the case $K(z)|z|^{-(n + \alpha)}$ has finite measure (see~\cite{Chasseigne}) to
the fractional Laplacian of order $\alpha$, which is the case when $K$ is equal to a  well-known
constant $C_{n,\alpha} > 0$ (see~\cite{Hitch}). 


Our main interest is to prove the well-posedness of problem~\eqref{cauchy} in the context of loss of the boundary condition, namely existence and 
uniqueness of a viscosity solution in $C(\bar{Q})$ which does not agree with $\varphi$ on $\partial \Omega \times (0,+\infty)$. 
Such losses of boundary conditions were studied in \cite{Barles-Chasseigne-Imbert} whose main result was that, if $H$ has some natural 
growth depending on the ellipticity properties of $\mathcal{I}$, then there is no loss of boundary condition. Our key assumptions on $\alpha$
and $H$ will imply that our framework is exactly the opposite, i.e. the $H$ term will be (in a suitable sense) stronger than the $\mathcal{I}$ one.

We recall that, in the second-order case, there are two well-known examples of problems developing this kind of loss of boundary conditions. 
The first case is the case of the degenerate parabolic problems where the $\mathcal{I}$ is replaced by a second-order linear operator: the equation becomes
\begin{equation*}
\partial_t u - \frac{1}{2} {\rm Tr}(a(x)D^2u(x)) + H(x,t,u,Du) = 0 \quad \mbox{in} \ \Omega \times [0,+\infty),
\end{equation*}
but we assume that the operator is degenerate, i.e. the symmetric matrix $a(x)$ is nonnegative for any $x$ but can have $0$ eigenvalues. 
Such problems, in particular in the linear case where studied by  Keldysh~\cite{Keldysh} and Radkevich~\cite{Radkevich1, Radkevich2} by pde methods (solutions in a weak sense) and by Freidlin~\cite{Freidlin} through a probabilistic approach. 
The first general results by a viscosity solutions' approach handling real losses of Dirichlet boundary conditions for second-order equations 
appears in \cite{Barles-Burdeau} following some previous results for first-order equations (see \cite{BP1,BP2}). More specifically, in problems which arise from the study of optimal exit time problems, one is led to Hamilton-Jacobi equations where $H$ has the Bellman form
\begin{equation}\label{Hamilton1}
H(x,t,u,p) = \lambda u + \sup \limits_{\beta \in \mathcal{B}} \{ -b(x,t, \beta) \cdot p - f_\beta(x,t, \beta)\}, 
\end{equation}
where  $\lambda \geq 0$, $\mathcal{B}$ is a compact metric space (the control-space) and $b, f$ 
are continuous and bounded functions (see~\cite{Bardi-Capuzzo} and \cite{Fleming-Soner} for the connections between control problems and such equations).

Loss of boundary conditions may arise at some point $x_0 \in \bar{\Omega}$ when $a(x_0)$ is singular, and more precisely when $a(x_0) n(x_0) = 0$ 
where $n(x_0)$ the unit outer normal vector to $\partial \Omega$ at $x_0$. This condition indicates the lack of diffusion in the normal 
direction at $x_0$. In this context, in order to decide if there is (or not) a loss of boundary condition, one has to examine the first-order term in 
the equation together with the geometrical properties of the boundary : we do not give details here and refer instead to \cite{Barles-Burdeau}.
Despite of the difficulty connected to the loss of boundary conditions, existence and uniqueness for such problems
can be obtained in the context of viscosity solutions with \textsl{generalized boundary condition} 
(see~\cite{Barles-Burdeau},~\cite{Barles-DaLio},~\cite{DaLio}~\cite{Barles-Rouy} and references therein). 

The second example, and in some sense which can be seen as being closer to our framework,  is the case of uniformly parabolic second-order problem associated to a Hamiltonian with 
superquadratic growth in $Du$, namely equations with the form
\begin{equation}\label{eqsuperquad}
\partial_t u - \Delta u + H(x,u,Du) = 0 \quad \mbox{in} \ \Omega \times [0,+\infty),
\end{equation}
where 
\begin{equation}\label{superquad}
H(x,t,u,p) = \lambda u + |p|^m - f(x), \quad m > 2,
\end{equation}
where $\lambda \geq 0$ and $f \in C(\bar{\Omega})$.
In this case, losses of boundary conditions come from the relative strenght of the second-order term and the $|Du|^m$-term : in the superquadratic case,
the $|Du|^m$-term may impose such losses of boundary data. In~\cite{Barles-DaLio_On},~\cite{Tchamba}, the existence and uniqueness of solutions is obtained (taking into account these losses of Dirichlet boundary conditions) and the asymptotic behavior of the solution of the problem as $t \to +\infty$ is also studied in \cite{Tchamba}. In this task, the discount rate $\lambda$ in problems with Hamiltonians as~\eqref{Hamilton1} or~\eqref{superquad} is determinant on the asymptotic behavior. For instance, as it can be seen in~\cite{Tchamba}, if $\lambda > 0$ then the asymptotic behavior of problems 
like~\eqref{eqsuperquad} is the uniform convergence in $C(\bar{\Omega})$ as $t \to +\infty$ to the solutions of the associated stationary problem.
However, if the case $\lambda = 0$ different behaviors may arise and it is well-known that the \textsl{ergodic problem} plays a key role, 
see~\cite{Barles-Souganidis}. 
We mention here that such as ergodic behavior for nonlocal operators is studied by the 
authors in collaboration with S. Koike and O. Ley~\cite{BKLT}, see also~\cite{Barles-Chasseigne-Ciomaga-Imbert}. 

This (very brief and incomplete) state-of-the-art on parabolic Dirichlet problems with loss of boundary conditions allows us to be 
more specific on the contents of this paper : we obtain the well-posedness of problem (CP) in two cases which can be understood as 
the extension of the both types of second-order problems we presented above. The first one concerns {\em coercive} Hamiltonians 
as~\eqref{superquad} for which the superquadratic condition has to be replaced in our context by the \textsl{superfractional condition}
$m > \alpha$, making the first-order term the leading term in the equation. We remark that we have no other additional restriction to $m$ 
(in particular, we can deal with $m<1$) and then we allow the study of Hamiltonians which are concave in $Du$.

On the other hand, in the case of problem~\eqref{cauchy} associated to Bellman-type Hamiltonians with the form~\eqref{Hamilton1}, 
the diffusive role of $\mathcal{I}$ defined in~\eqref{operator} is of weaker order than the first-order term when we assume 
$\alpha < 1$. However, in contrast with the {\em degenerate} second-order case, losses of boundary conditions arise even if we impose an
\textsl{uniform ellipticity condition} in the sense of Caffarelli 
and Silvestre~\cite{Caffarelli-Silvestre}, which is related with the nonintegrability of $K^\alpha$ at the origin (see assumption (UE) below).
As in \cite{Topp}, the well-posedness of (CP) is obtained 
through a careful examination of the effects of the drift $b$ at each point of $\partial \Omega \times (0,+\infty)$ and suitable assumptions.


\medskip
\noindent
{\bf Organization of the Paper:} In Section~\ref{defsolsection} we provide the notion of solution for~\eqref{cauchy}. 
In section~\ref{assumptionsection} we precise what we mean with~\eqref{cauchy} in coercive and Bellman form,
introduce the assumptions of each problem and present the main results. In section~\ref{bcsection} we study the behavior of sub and 
supersolutions on the parabolic boundary. Section~\ref{regularitysection} is devoted to regularity issues for each problem.
The proof of the main results are given in section~\ref{proofsection} and the existence, uniqueness and large time behavior 
is addressed in section~\ref{existencesection}.


\section{Basic Notation and Notion of Solution.}
\label{defsolsection}

We start with the basic notation.
For $\delta > 0$ and $x \in \R^n$ we write $B_\delta(x)$ as the ball of radius $\delta$ centered at $x$ and $B_\delta$ if $x = 0$.
For an arbitrary set $A$, we denote $d_A(x) = dist(x, \partial A)$ the signed distance function to $\partial A$ which is nonnegative for 
$x \in A$ and nonpositive for $x \notin A$. For $\Omega$ we simply
write $d(x) = d_{\partial \Omega}(x)$ and define the set $\Omega_\delta$ as the open set of all 
$x \in \Omega$ such that $d(x) < \delta$. By the smoothness of the domain, there exists a fixed number $\delta_0 > 0$,
depending only on $\Omega$, such that $d$ is smooth in the set of points $x$ such that $|d(x)| < \delta_0$ (see~\cite{G-T}). 
For $x \in \R^n$ and $\lambda \in \R$, we write
$$
\Omega - x = \{ z  : x + z \in \Omega\} \quad \mbox{and} \quad \lambda \Omega = \{ \lambda z : z \in \Omega\}. 
$$

By a modulus of continuity $\omega$ we mean a nondecreasing, sublinear, continuous function $\omega: [0,+\infty) \to \R$ such that $\omega(0) = 0$.

Given a set $A \subset \R^n$, we denote $\mathrm{USC}(A)$ the set of real valued, upper semicontinuous (usc for short) functions. In the analogous
way, we write $\mathrm{LSC}(A)$ the set of real valued, lower semicontinuous (lsc for short) functions.

Before presenting the viscosity evaluation, we need to introduce some notation related with the nonlocal term $\I$.
For $\alpha \in (0,2)$, we denote 
$$
K^\alpha(z) = K(z)|z|^{-(n + \alpha)}, \quad \mbox{for } \ z \neq 0.
$$ 

As we mentioned in the introduction, we are interested in the case $\alpha$ represents the order of $\I$ and therefore, in the case 
$\alpha \in (0,1)$, for each $x \in \R^n$ and $\phi: \R^n \to \R$ bounded and smooth at $x$, we write
\begin{equation}\label{operatoralphasmall}
\I(\phi, x) = \int_{\R^n} [\phi(x + z) - \phi(x)] K^\alpha(z)dz.
\end{equation}

We remark that in the case $K$ is symmetric (that is, $K(z) = K(-z)$ for $z \in \R^n$), then~\eqref{operator}
is equivalent to~\eqref{operatoralphasmall} when $\alpha \in (0,1)$.

For $x, p \in \R^n$, $A \subset \R^n$ and $\phi$ a bounded function, we define
\begin{equation}\label{defIdelta}
\I[A](\phi, x, p) = \int_{\R^n \cap A} [\phi(x + z) - \phi(x) - \mathbf{1}_{B} \langle p, z \rangle] K^\alpha(z)dz.
\end{equation}

We write in a simpler way $\I[A](\phi, x) = \I[A](\phi, x, D\phi(x))$ when $\phi \in L^\infty(\R^n) \cap C^2(B_\delta)$ for some $\delta > 0$,
$\I(\phi, x, p) = \I[\R^n](\phi, x, p)$ when $A = \R^n$. 
In the case $\alpha \in (0,1)$, the presence of the compensator (namely, the term $\mathbf{1}_B \langle D\phi(x), z \rangle$)
is not necessary to give a sense to the nonlocal term and for this reason we drop it in~\eqref{defIdelta}.
 
If $\phi \in C^2(B_{\delta}(x) \times (t - \delta, t +\delta))$ and 
$w : \R^n \times \R \to \R$ is a bounded measurable function, we define
\begin{equation}\label{Eevaluation}
\begin{split}
& E_{\delta}(w, \phi, x, t) \\
= & \ \partial_t \phi(x,t) - \I[B_\delta](\phi(\cdot, t), x) - \I[B_\delta^c](w(\cdot, t), x, D\phi(x,t)) \\
& + H(x,t, w(x,t), D\phi(x,t)),
\end{split}
\end{equation}
where ``$E$'' stands for ``evaluation''. 

For $T > 0$, we define the sets
\begin{equation*}
Q_T = \Omega \times (0,T]; \quad \pboundaryT = \partial \Omega \times (0,T]; \quad \pexteriorT = \Omega^c \times (0,T].  
\end{equation*}

We are going to consider finite time horizon problem associated with~\eqref{cauchy}
\begin{equation}\tag{C$\rm{P}_T$}\label{cauchyT}
\left \{ \begin{array}{rll} \partial_t u - \mathcal{I}[u(\cdot, t)] + H(x, t, u, D u) &= 0 \ &\mbox{in} \ Q_T \\
u(x,t) &= \varphi(x,t) \ &\mbox{in} \ \pexteriorT \\
u(x,0) &= u_0(x) \ &\mbox{in} \ \bar{\Omega}, \end{array} \right .
\end{equation}

For a function $u \in \mathrm{USC}(\bar{\Omega} \times [0,T])$ (resp $u \in \mathrm{LSC}(\bar{\Omega} \times [0,T])$), 
we define its upper (resp. lower) $\varphi$-extension as the function defined in $\R^n \times [0,T]$ by
\begin{equation}\label{Uextension}
\begin{split}
& u^\varphi(x,t) \ (\mbox{resp.} \ u_\varphi(x, t)) \\
= & 
\left \{ \begin{array}{ll} u(x, t) & \mbox{if} \ (x,t) \in \Omega \times [0,T] \\ \varphi(x, t) & \mbox{if} \ (x,t) \in \bar{\Omega}^c \times [0,T] \\ 
\max \ (\mbox{resp.} \ \min) \{ u(x, t), \varphi(x, t) \} & \mbox{if} \ (x,t) \in \partial \Omega \times [0,T], \end{array} \right . 
\end{split}
\end{equation}

We provide a  definition of solution to problem~\eqref{cauchyT}
which can be extended naturally to~\eqref{cauchy}. 
\begin{defi}\label{defsol}
A function $u \in \mathrm{USC}(\bar{\Omega} \times [0,T])$ is a viscosity subsolution 
of~\eqref{cauchyT} if for any smooth function $\phi: \R^n \times [0,T] \to \R$, any maximum point 
$(x_0, t_0) \in \bar{\Omega} \times [0,T]$ of $u^\varphi - \phi$ in $B_\delta(x_0) \times (t_0 - \delta, t_0 + \delta) \cap \R^n \times [0,T]$
with $\delta > 0$, we have the inequality
\begin{eqnarray*}
E_{\delta}(u^\varphi, \phi, x_0, t_0) \leq 0 && \quad \mbox{if} \  (x_0, t_0) \in Q_T, \\
\min \{ E_{\delta}(u^\varphi, \phi, x_0, t_0), u(x_0, t_0) - \varphi(x_0, t_0) \} \leq 0 && \quad \mbox{if} \ x_0 \in \partial \Omega, \\
\min \{ E_{\delta}(u^\varphi, \phi, x_0, t_0), u(x_0, t_0) - u_0(x_0) \} \leq 0 && \quad \mbox{if} \ t_0 = 0,
\end{eqnarray*}
where $E_\delta$ is defined in~\eqref{Eevaluation}.

A function $v \in \mathrm{LSC}(\bar{\Omega} \times [0,T])$ is a viscosity supersolution 
of~\eqref{cauchyT} if for any smooth function $\phi: \R^n \times [0,T] \to \R$, any minimum point 
$(x_0, t_0) \in \bar{\Omega} \times [0,T]$ of $v_\varphi - \phi$ in $B_\delta(x_0) \times (t_0 - \delta, t_0 + \delta) \cap \R^n \times [0,T]$
with $\delta > 0$, we have the inequality
\begin{eqnarray*}
E_{\delta}(v_\varphi, \phi, x_0, t_0) \geq 0 && \quad \mbox{if} \  (x_0, t_0) \in Q_T, \\
\max \{ E_{\delta}(v_\varphi, \phi, x_0, t_0), v(x_0, t_0) - \varphi(x_0, t_0) \} \geq 0 && \quad \mbox{if} \ x_0 \in \partial \Omega, \\
\max \{ E_{\delta}(u^\varphi, \phi, x_0, t_0), u(x_0, t_0) - u_0(x_0) \} \geq 0 && \quad \mbox{if} \ t_0 = 0.
\end{eqnarray*}

Finally, a viscosity solution of~\eqref{cauchyT} is a function whose upper and lower semicontinuous envelopes are sub and supersolution 
of the problem, respectively. 
\end{defi}

The above definition is basically the same as the one presented in~\cite{Alvarez-Tourin}, \cite{Barles-Chasseigne-Imbert}, \cite{Barles-Imbert},~\cite{Sayah1} 
and~\cite{Sayah2}. Written in that way we highlight the goal of this paper, which is to state the existence and uniqueness of a solution 
of~\eqref{cauchy} in $C(\bar{Q})$. 

We note that Definition~\ref{defsol} interprets the points at $\Omega \times \{ T \}$ as \textsl{interior points}, which is consistent with 
the classical definition of the Cauchy problem for parabolic equations (see~\cite{Evans},\cite{Friedman}). Of course, a weaker definition 
of viscosity solution (concerning functions defined only in $\bar{\Omega} \times [0,T)$) can be set, obtaining the same results 
presented in this paper. However, we avoid this extra difficulty here since its consideration has no significant contribution
to the development of our problem.


\section{Assumptions and Main Results.}
\label{assumptionsection}

As we mentioned in the introduction, in this paper we study the well-posedness for problem~\eqref{cauchy} in two cases, depending on the
features of $H$. Basically, we are interested in the case when $H$ has a coercive nature in the gradient term, 
and the case $H$ has a Bellman form and therefore it is not necessarily coercive.

\subsection{Coercive Hamiltonian and Examples.}
In this case we restrict the time dependence of $H$ by the assumption

\medskip
\noindent
(A0) \textsl{ There exists $H_0: \bar{\Omega} \times \R \times \R^n \to \R$ continuous and $f : \bar{\Omega} \times [0,+\infty) \to \R$ 
uniformly continuous and bounded such that
\begin{equation*}
H(x,t,r,p) = H_0(x,r,p) - f(x,t), 
\end{equation*}
for all $x \in \bar{\Omega}, t \geq 0, r \in \R$ and $p \in \R^n$.}

\medskip

Let $\alpha \in (0,2)$ and $\I$ as in~\eqref{operator},~\eqref{operatoralphasmall}. We will consider \textsl{superfractional}
coercive Hamiltonians, where the gradient growth is given by $H_0$ through the basic assumption

\medskip

\noindent
(A1) \textsl{ There exists $m > \alpha$ and $C_0 > 0$ such that, for all $R > 0$ there exists $C_R > 0$ satisfying
\begin{equation*}
H_0(x,r,p) \geq C_0 |p|^m - C_R, 
\end{equation*}
for all  $x \in \bar{\Omega}$, $p \in \R^n$ and $|r| \leq R$.}

\medskip

However, we must be careful if the coercivity is sub or superlinear. For this, we split the analysis depending on the gradient growth of $H_0$, that is

\medskip
\noindent
$\bullet$ \textsl{Sublinear Coercivity:} Assume (A0) holds. We say that $H$ is \textsl{sublinearly coercive} if 
$H_0$ satisfies (A1) with $m \leq 1$, and the following continuity condition holds

\medskip
\noindent
(A2-a) \textsl{For all $R > 0$, there exists a modulus of continuity $\omega_R$ satisfying}
\begin{equation*}
\begin{split}
H_0(y, r, p) - H_0(x, r, p + q)
\leq  \omega_R( |x - y|(1 +  |p|)  + |q|),
\end{split}
\end{equation*}
\textsl{for all $x, y \in \bar{\Omega}$, $|r| \leq R$, $p,q \in \R^n$, $|q|\leq 1$.}

\medskip
\noindent
$\bullet$ \textsl{Superlinear Coercivity:} Assume (A0) holds. We say that $H$ is \textsl{superlinearly coercive} if 
$H_0$ satisfies

\medskip

\noindent
(A1-b) \textsl{ There exists $m > \max \{1, \alpha \}$ and $a_0 > 0$ such that, for all $R > 0$, there exists a constant $C_R$ such that
\begin{equation*}
H_0(x,r, p) - \mu H_0(x, \mu^{-1}r, \mu^{-1} p) \leq \Big{(} (1 - m) a_0 |p|^m + C_R \Big{)} (1 - \mu)
\end{equation*}
for all $\mu < 1$, $x \in \bar{\Omega}$, $p \in \R^n$, $|r| \leq R$.}

\medskip
\noindent
(A2-b) \textsl{If $m$ is given by Assumption (A1-b), for all $R > 0$, there exists a modulus of continuity $\omega_R$ satisfying}
\begin{equation*}
\begin{split}
H_0(y, r, p) - H_0(x, r, p + q)
\leq  \omega_R(|x - y|)(1 +  |p|^m) + |p|^{m - 1} \omega_R(|q|),
\end{split}
\end{equation*}
\textsl{for all $x, y \in \bar{\Omega}$, $|r| \leq R$, $p,q \in \R^n$ , $|q|\leq 1$.}

\begin{remark}\label{rmksupercoercive}
Note that Condition (A1-b) gives us the gradient coercivity of $H_0$ since it implies (A1) with $m > \max \{ 1, \alpha \}$. 
\end{remark}

In order to describe the kind of Hamiltonians we have in mind, we introduce the following examples : in the first one, we assume $m\leq 1$ and consider
\begin{equation*}
H(x,t, r,p) = a_1(x) |p|^m + a_2(x)|p|^l + \lambda(x) r - f(x, t),
\end{equation*}
while in the second case, we suppose $m > 1$ and
\begin{equation*}
H(x,t, r,p) = a_1(x) |p|^m + a_2(x)|p|^l + b(x) \cdot p + \lambda(x) r - f(x, t).
\end{equation*}

In both cases, $l < m$, $a_1, a_2, \lambda, f: \bar{\Omega} \to \R$ are continuous functions 
with $\lambda \geq 0$. We assume in addition that $a_1, a_2$ are Lipschitz continuous and $a_1 \geq C_0$ for some fixed constant $C_0 > 0$. 

These Hamiltonians are coercive in $Du$ and in the case $m>1$ we can include transport terms
with a Lipschitz continuous vector field $b: \bar{\Omega} \to \R^n$. The above assumptions are easily checkable in both cases.

\subsection{Bellman Hamiltonian.}
Let $\mathcal{B}$ a compact metric space, 
$b: \bar{\Omega} \times [0,+\infty) \times \mathcal{B} \to \R^n$ and $f, \lambda: \bar{\Omega} \times [0,+\infty) \times \mathcal{B} \to \R$ 
continuous and bounded functions. We say that $H$ has a \textsl{Bellman form} if, for $t \in [0,+\infty)$, $r \in \R$, $x \in \bar{\Omega},p \in \R^n$, 
$H(x,t,r,p)$ can be written as
\begin{equation}\tag{$H_\mathcal{B}$}\label{Hcontrol}
H(x,t,r, p) = \sup \limits_{\beta \in \mathcal{B}} \{\lambda_\beta(x,t) r - b_\beta(x,t) \cdot p - f_\beta(x,t) \}, 
\end{equation}
and satisfies the assumptions~\eqref{L} and~\eqref{Gamma} below. In~\eqref{Hcontrol} 
we have adopted the abuse of notation $b_\beta(x,t) = b(x, t, \beta)$  and in the same way for the other functions. 

For $H$ with the form~\eqref{Hcontrol} we impose the uniform space-time Lipschitz asumption:

\medskip

\noindent
\textsl{There exists $L > 0$ such that, for all $\beta \in \mathcal{B}$, $(x,s), (y,t) \in \bar{\Omega} \times [0,+\infty))$, we have}
\begin{equation}\tag{L}\label{L}
 |b_\beta(x,s) - b_\beta(y,t)| \leq L(|x - y| + |s - t|).
\end{equation}

\medskip

Then we introduce the notation
\begin{eqnarray*}
\begin{split}
\Gamma_{in} & = \{ (x,t) \in \pboundary \ : \ \forall \ \beta \in \mathcal{B}, \ b_\beta(x,t) \cdot Dd(x) > 0 \}, \\ 
\Gamma_{out} & = \{ (x,t) \in \pboundary \ : \ \forall \beta \in \mathcal{B}, \ b_\beta(x,t) \cdot Dd(x) \leq 0 \}, \\
\Gamma & = \pboundary \setminus (\Gamma_{in} \cup \Gamma_{out}),
\end{split}
\end{eqnarray*}
and with this, we consider the following condition over the behavior of the drisf terms on $\pboundary$
\textsl{
\begin{equation}\tag{$\Sigma$} \label{Gamma}
\Gamma_{in}, \Gamma_{out} \ \mbox{and} \ \Gamma \ \mbox{are unions of connected components of} \ \pboundary. 
\end{equation}
}

We remark that, in the current Bellman setting, the nonlocal term $\I$ is assumed to be of order $\alpha < 1$. Therefore it has a weaker effect compared with the first-order terms. In particular, on the boundary, the behavior of the drift  plays a determinant role.
In this direction, the set $\Gamma_{out}$ should be understood as the set where the classical boundary condition holds, 
meanwhile on $\Gamma_{in}$ may arise losses of the boundary condition due to the ``stronger'' influence of the transport term compared with the nonlocal diffusion. 
Finally, on $\Gamma$, we do not have a transport effect anymore : the value of the different costs (boundary or running cost) decides of the choice of the control and of the loss or no loss of boundary condition.

We introduce assumption~\eqref{Gamma} in order to avoid have different behaviors of the $b_\beta$'s on the 
same connected component, which could be a source of discontinuities for the solution 
(the reader may think in term of transport equation to be convinced by this claim). On $\Gamma$, it can be seen as a controllability assumption in the normal direction.
Similar assumptions of the boundary are made in~\cite{Barles-Burdeau},~\cite{Barles-Rouy} in the degenerate second-order setting and~\cite{Topp} for 
the nonlocal one.

\subsection{Structural Assumptions and Main Results.}

As it is classical for Cauchy-Dirichlet problems, the initial and boundary data satisfy the following \textsl{compatibility condition} 
at $t = 0$

\medskip
\noindent
(H0) \textsl{ \qquad \qquad \qquad \qquad $u_0(x) = \varphi(x, 0)$, for all  $x \in \partial \Omega$.
}

\medskip

The properness of the problem is encoded by the following two conditions

\medskip
\noindent
(H1) \textsl{For all $R > 0$, there exists $h_R \in C(\bar{\Omega})$ such that, for all 
$x \in \bar{\Omega}$, $u,v \in \R$, $0 \leq t \leq R$, and $p \in \R^n$, we have
\begin{equation*}
H(x,t,u,p) - H(x,t,v,p) \geq h_R(x) (u - v).
\end{equation*}
}

\medskip
\noindent
(H2) \textsl{For all $R > 0$, the function $h_R$ in (H1) satisfies
\begin{equation*}
\inf \limits_{x \in \bar{\Omega}} \{ h_R(x) + \int \limits_{\Omega^c - x} K^\alpha(z)dz \} \geq 0.
\end{equation*}
}

As it is classical in problems where loss of the boundary condition arises, 
Strong Comparison Principle needs the introduction of a modification of sub and 
supersolutions. For a function $u$ bounded
and usc in $\bar{Q}$ (which will be thought as subsolution) we denote
\begin{equation}\label{deftildeuvm}
\begin{split}
\tilde{u}(x, t)  = 
\left \{ \begin{array}{lccl} &\limsup \limits_{\begin{tiny} Q \ni (y,s) \to (x,t) \end{tiny}} 
u(y, s)& \quad & \mbox{if} \ (x, t) \in \pboundary \\
&u(x, t)& \quad & \mbox{if} \ (x, t) \in \bar{Q} \setminus \pboundary.  \end{array} \right .
\end{split}
\end{equation}

\begin{teo}\label{teocoercive}{\bf (Strong Comparison Principle - Coercive Case)}
Let $\varphi \in C_b(\bar{Q}^{ext})$ and $u_0 \in C(\bar{\Omega})$. Assume (H0) holds and that $H$ 
has a coercive form satisfying (H1)-(H2). 
If $u, v$ are bounded viscosity sub and supersolution to problem~\eqref{cauchy} respectively, then
\begin{equation*}
u \leq v \quad \mbox{in} \ Q \cup \bar{\Omega} \times \{ 0 \}. 
\end{equation*}

Moreover, if $\tilde{u}$ is defined as in~\eqref{deftildeuvm}, then $\tilde{u} \leq v$ in $\bar{Q}$.
\end{teo}

The result concerning the Bellman needs also a redefinition of sub and supersolutions at the boundary.
Of course, in this control framework, the different part of the boundary $\{ \Gamma_{in}, \Gamma_{out}, \Gamma \}$ play different roles.

For bounded functions $u$ and $v$, $u$ usc in $\bar{Q}$, $v$ lsc in $\bar{Q}$, we denote
\begin{equation}\label{deftildeuv}
\begin{split}
\tilde{u}(x, t) &= 
\left \{ \begin{array}{lccl} &u(x, t)& \quad & \mbox{if} \ (x, t) \in \bar{Q} \setminus (\Gamma_{in} \cup \Gamma) \\
&\limsup \limits_{\begin{tiny} Q \ni (y,s) \to (x,t) \end{tiny}} 
u(y, s)& \quad & \mbox{if} \ (x, t) \in \Gamma_{in} \cup \Gamma  \end{array} \right . \\
\tilde{v}(x, t) &= 
\left \{ \begin{array}{lccl} &v(x,t)& \quad & \mbox{if} \ (x,t) \in \bar{Q} \setminus \Gamma_{in} \\
&\liminf \limits_{\begin{tiny} Q \ni (y,s) \to (x,t) \end{tiny}} 
v(y, s)& \quad & \mbox{if} \ (x,t) \in \Gamma_{in}.  \end{array} \right .
\end{split}
\end{equation}

In the Bellman case, we will require the stronger ellipticity assumption
\begin{equation}\tag{UE}\label{uellipticity}
\begin{split}
\mbox{There exists} \ c_1, c_2 > 0 \ \mbox{such that} \ c_2 \leq K(z) \ \mbox{for all} \ |z| \leq c_1.
\end{split}
\end{equation} 

\begin{teo}\label{teoBellman}{\bf (Strong Comparison Principle - Bellman Case)}
Let $\varphi \in C_b(\bar{Q}^{ext})$ and $u_0 \in C(\bar{\Omega})$. Assume $\alpha < 1$, (UE), (H0) hold and let $H$ with Bellman form
satisfying (H1)-(H2).  
If $u, v$ are bounded viscosity sub and supersolution of \eqref{cauchy} respectively, then
\begin{equation*}
u \leq v \quad \mbox{in} \ Q \cup \bar{\Omega} \times \{ 0 \}. 
\end{equation*}

Moreover, if $\tilde{u}, \tilde{v}$ are defined as in~\eqref{deftildeuv}, then $\tilde{u} \leq \tilde{v}$ in $\bar{Q}$.
\end{teo}

The result of Theorem~\ref{teoBellman} can be obtained without the uniform ellipticity assumption (UE) by slightly changing the 
definition of $\Gamma_{in}, \Gamma_{out}$ and $\Gamma$. Indeed, in this setting, only the assumptions on the drift term determine the 
loss or not loss of the boundary condition of the solution on $\Gamma_{in}, \Gamma_{out}$ and $\Gamma$ and they have to be strong enough to 
compensate the lack of the ellipticity effect of $\I$.


\section{Initial and Boundary Condition.}
\label{bcsection}

We also remark that, considered as a part of the parabolic boundary, we ask the initial condition is satisfied in the generalized sense. However, 
the initial condition is satisfied in the classical sense on $\Omega \times \{ 0 \}$. Moreover,
mainly because of~(H0), the condition holds classically on $\bar{\Omega} \times \{ 0 \}$.
\begin{lema}\label{lemat=0}
Assume that $H \in C(\bar{\Omega} \times [0,+\infty) \times \R \times \R^n)$, $\varphi \in C_b(\bar{Q}^{ext})$, $u_0 \in C(\bar{\Omega})$
satisfying~(H0). 
If $u, v$ are respectively a bounded, usc viscosity subsolution and a bounded, lsc viscosity supersolution to~\eqref{cauchy}, then $u(x, 0) \leq u_0(x) \leq v(x, 0)$ for all $x \in \bar{\Omega}$.
\end{lema}

The proof of this lemma follows the same lines of the analogous result for the second-order case presented in~\cite{DaLio}, with subtle 
modifications concerning the nonlocal operator.

Now we look for the behavior of sub and supersolutions at the lateral parabolic boundary.
\begin{lema}\label{lateralboundary} Assume that $H \in C(\bar{\Omega} \times [0,+\infty) \times \R \times \R^n)$ and
$\varphi \in C_b(\bar{Q}^{ext})$. If $(x_0, t_0) \in \pboundary$ and $u, v$ are respectively a bounded, usc viscosity subsolution 
and a bounded, lsc viscosity supersolution to~\eqref{cauchy}, then

\noindent
{\bf (i)} We have $u(x_0, t_0) \leq \varphi(x_0, t_0)$ if one of the following conditions hold:

\begin{itemize}
 \item[(i.1)] There exists $C_0, \rho > 0$ and $m > \alpha$ such that for all $R > 0$, there exists $C_R > 0$ satisfying 
$$
H(x,t,r, k\eta^{-1} Dd(x)(1 + o_\eta(1))) \geq C_0 (k\eta^{-1})^m - C_R
$$ 
for all $k, \eta > 0$, $x \in B_\rho(x_0)$ and $t, |r| \leq R$. 

 \item[(i.2)] Condition~\eqref{uellipticity} with $\alpha < 1$ holds, 
and there exists $c_0, \rho > 0$ such that, for all $R > 0$ there exists $C_R$ satisfying 
$$
H(x, t, r, k\eta^{-1} Dd(x)(1 + o_\eta(1))) \geq -c_0 k \eta^{-1} d(x) - C_R
$$ 
for all $k, \eta > 0$, $x \in B_\rho(x_0)$ and $t, |r| \leq R$.
\end{itemize}

\noindent
{\bf (ii)} We have $v(x_0, t_0) \geq \varphi(x_0, t_0)$ if condition~\eqref{uellipticity} with $\alpha < 1$ holds,
and there exists $c_0, \rho > 0$ such that, for all $R > 0$ there exists $C_R$ satisfying 
$$
H(x, t, r, -k\eta^{-1} Dd(x)(1 + o_\eta(1))) \leq c_0 k\eta^{-1} d(x) + C_R
$$ 
for all $k, \eta > 0$, $x \in B_\rho(x_0)$ and $t, |r| \leq R$.
\end{lema}


\noindent
{\bf \textit{Proof:}} We concentrate on {\bf (i)} since {\bf (ii)} is an adaptation to $(i.2)$. 
By contradiction, we assume $u(x_0, t_0) - \varphi(x_0, t_0) = \nu$ for some $\nu > 0$. This  
implies in particular that $u^\varphi(x_0, t_0) = u(x_0, t_0)$. We consider $\sigma \in (\max\{ 1, \alpha \}, 2)$ and $C^{1, \sigma - 1}$ functions
$\chi, \psi: \R \to \R$ such 
that $\chi$ is even, bounded, $\chi(0) = 0$, $\chi(t) > 0$ for $t \neq 0$, $\liminf_{|t| \to \infty} \chi(t) > 0$ and such
that $\chi(t) = |t|^\sigma$ in a neighborhood of $0$. For $\psi$ we assume it is
bounded, strictly increasing, $\psi \geq -\frac{1}{4} \nu$ and such that for some $k > 0$, $\psi(t) = k t$ for all $|t| \leq 1$. We consider a 
parameter $\eta$ and $\epsilon = \epsilon_\eta \to 0$ as $\eta \to 0$ to be fixed later, and introduce the test function
$$
\Psi(y, t) := \psi(d(y)/\eta) + \epsilon^{-1} \chi(|y - x_0|) + \epsilon^{-1} \chi(|t - t_0|).
$$

By our assumption on $u, \varphi$, $\chi$ and $\psi$, the function $(x,t) \mapsto u^\varphi(x,t) - \Psi(x,t)$ has a maximum point 
$(\bar{x}, \bar{t}) \in \R^n \times (0,T)$ for $\eta$ small enough. Of course, $(\bar{x}, \bar{t})$ depends on $\eta$ but we drop the dependence on 
$\eta$ to simplify the notations. From the maximum point property, $u^\varphi(\bar{x}, \bar{t}) - \Psi(\bar{x}, \bar{t}) \geq u^\varphi(x_0, t_0) - \Psi(x_0, t_0)$ which implies
\begin{equation*}
u^\varphi(\bar{x}, \bar{t}) - \psi(d(\bar{x})/\eta) 
- \epsilon^{-1} \chi(|\bar{x} - x_0|) - \epsilon^{-1} \chi(|\bar{t} - t_0|) \geq \varphi(x_0) + \nu.
\end{equation*}
Using this inequality, classical arguments show that $\bar{x} \to x_0$ and $\bar{t} \to t_0$ as $\eta \to 0$. And from the same inequality we obtain 
$\bar{x} \in \bar{\Omega}$ for $\eta$ small enough because $\psi \geq -1/4 \nu$ and $\varphi$ is continuous. Finally, using properly the usc of $u^\varphi$ we conclude
\begin{equation}\label{proplateral}
d(\bar{x}) = o_1(\eta) \eta, \ |\bar{x} - x_0|, |\bar{t} - t_0| = o_\eta(1), \  \mbox{and} \ 
u^\varphi(\bar{x}, \bar{t}) \to u(x_0, t_0),
\end{equation}
as $\eta \to 0$. Hence, picking some $\delta >0$, we can use the viscosity inequality for subsolutions, concluding that
\begin{equation}\label{testinglateral}
\begin{split}
\partial_t \Psi(\bar{x}, \bar{t}) \leq & \ \mathcal{I}[B_\delta](\Psi(\cdot, \bar{t}), \bar{x}) + 
\mathcal{I}[B_\delta^c](u^\varphi(\cdot, \bar{t}), \bar{x}, D\Psi(\bar{x}, \bar{t})) \\
& \ - H(\bar{x}, \bar{t}, u(\bar{x}, \bar{t}), D\Psi(\bar{x}, \bar{t})), 
\end{split}
\end{equation}
where in view of the first and second statement in~\eqref{proplateral}, for $\eta$ small enough we can write
\begin{equation}\label{DPsi}
D\Psi(\bar{x}, \bar{t}) = k \eta^{-1} Dd(\bar{x}) + \epsilon^{-1} |\bar{x} - x_0|^{\sigma - 2} (\bar{x} - x_0).
\end{equation}
 
We start with the estimates concerning the nonlocal terms in~\eqref{testinglateral}. 
To do this, we consider $r \leq 1$ independent of $\eta$ and $d(\bar{x}) < \delta \leq \mu < r$. We define the sets
\begin{equation*}
\begin{split}
& \mathcal{A}_{\delta}^{ext} = \{ z \in B_r \ : \ d(\bar{x} + z) \leq d(\bar{x}) - \delta \}. \\ 
& \mathcal{A}_{\delta, \mu} = \{ z \in B_r \ : \ d(\bar{x}) - \delta < d(\bar{x} + z) < d(\bar{x}) + \mu \}. \\ 
& \mathcal{A}_{\mu}^{int} = \{ z \in B_r \ : \ \mu + d(\bar{x}) \leq d(\bar{x} + z) \}. 
\end{split}
\end{equation*}

We remark that $B_\delta \subset \mathcal{A}_{\delta, \mu}$ and 
using that $\bar{x}$ is a global maximum point of $u - \Psi$, in particular we
have $\delta(u^\varphi(\cdot, \bar{t}), \bar{x}, z) \leq \delta(\Psi(\cdot, \bar{t}), \bar{x}, z)$ in $\mathcal{A}_{\delta, \mu} \setminus B_\delta$.
Using this last fact we can write
\begin{equation*}
\begin{split}
& \I[B_\delta^c](u^\varphi(\cdot, \bar{t}), \bar{x}, D\Psi(\bar{x}, \bar{t})) + \mathcal{I}[B_\delta](\Psi(\cdot, \bar{t}), \bar{x}) \\
\leq & \ \I[B_r^c](u^\varphi(\cdot, \bar{t}), \bar{x}, D\Psi(\bar{x}, \bar{t})) 
+ \I[\mathcal{A}_{\delta}^{ext}](u^\varphi(\cdot, \bar{t}), \bar{x}, D\Psi(\bar{x}, \bar{t})) \\
& + \I[\mathcal{A}_{\mu}^{int}](u^\varphi(\cdot, \bar{t}), \bar{x}, D\Psi(\bar{x}, \bar{t})) 
+ \I[\mathcal{A}_{\delta, \mu}](\Psi(\cdot, \bar{t}), \bar{x}),
\end{split}
\end{equation*}
and from this we estimate each term in the right-hans side of the above inequality separately.
The constant $C > 0$ arising in each of the following estimates does not depend 
on $\mu, \delta, \eta$ or $\epsilon$. 

Using the expression~\eqref{DPsi}, we have
\begin{equation*}
\begin{split}
\I[B_r^c](u^\varphi(\cdot, \bar{t}), \bar{x}, D\Psi(\bar{x}, \bar{t})) \leq & \ 2||u^\varphi||_\infty \int_{B_r^c} K^\alpha(z)dz \\
& + (k\eta^{-1} + \epsilon^{-1} o_\eta(1)) \int_{B \setminus B_r} K^{\alpha - 1}(z),
\end{split}
\end{equation*}
where the last integral does not exists if $\alpha < 1$. Thus, we get
\begin{equation*}
\I[B_r^c](u^\varphi(\cdot, \bar{t}), \bar{x}, D\Psi(\bar{x}, \bar{t})) \leq C ||u^\varphi||_\infty r^{-\alpha} + 
C (\eta^{-1} + \epsilon^{-1} o_\eta(1)) r^{1 - \alpha},
\end{equation*}
and similarly, we have
\begin{equation*}
\I[\mathcal{A}_{\mu}^{int}](u^\varphi(\cdot, \bar{t}), \bar{x}, D\Psi(\bar{x}, \bar{t})) \leq C||u^\varphi||_\infty \mu^{-\alpha}
+ C (\eta^{-1} + \epsilon^{-1} o_\eta(1)) \mu^{1 - \alpha}.
\end{equation*}

At this point, we consider $\mu = \eta$. Thus, for all $\eta$ small enough and $z \in \mathcal{A}_{\delta, \mu}$ 
we have $\psi(d(\bar{x} + z)/\eta) = k\eta^{-1}d(\bar{x} + z)$ and applying the definition of $\Psi$ we get
\begin{equation*}
\begin{split}
& \Psi(\bar{x} + z, \bar{t}) - \Psi(\bar{x}, \bar{t}) \leq C(\eta^{-1} + \epsilon^{-1})|z|, \\
& \Psi(\bar{x} + z, \bar{t}) - \Psi(\bar{x}, \bar{t}) - \langle D\Psi(\bar{x}, \bar{t}), z \rangle \leq C(\eta^{-1} + \epsilon^{-1}) |z|^2,
\end{split}
\end{equation*}
from which we can get
\begin{equation*}\label{estimatelateral1}
\begin{split}
\I[\mathcal{A}_{\delta, \mu}](\Psi(\cdot, \bar{t}), \bar{x}) \leq C(\eta^{-1} + \epsilon^{-1}) \varrho_{\alpha}(\mu), 
\end{split}
\end{equation*}
where 
\begin{equation}
\varrho(\mu) = \left \{ \begin{array}{ll} \mu^{2 - \alpha} \quad & \mbox{if} \ \alpha > 1 \\ \mu \ln(\mu) \quad & \mbox{if} \ \alpha = 1
\\ \mu^{1 - \alpha} \quad & \mbox{if} \ \alpha < 1. \end{array} \right .
\end{equation}

Thus, recalling that we have chosen $\mu = \eta$ and taking $\epsilon \geq \eta^{\min \{\alpha, 1 \}}$, by the above estimates we can write
\begin{equation}\label{lateralnonlocal}
\begin{split}
& \I[B_\delta^c](u^\varphi(\cdot, \bar{t}), \bar{x}, D\Psi(\bar{x}, \bar{t})) + \mathcal{I}[B_\delta](\Psi(\cdot, \bar{t}), \bar{x}) \\
\leq & \ C \eta^{-\alpha} + C\eta^{-1} \varrho_{\alpha}(\mu) 
+ \I[\mathcal{A}_{\delta}^{ext}](u^\varphi(\cdot, \bar{t}), \bar{x}, D\Psi(\bar{x}, \bar{t})).
\end{split}
\end{equation}
where the constant $C$ depends only on the data and $||u^\varphi||_\infty$.

Under the above choice of $\epsilon$ and using~\eqref{proplateral}, we have 
$\partial_t \Psi(\bar{x}, \bar{t}) \geq \eta^{-\alpha} o_\eta(1)$. Using this estimate and~\eqref{lateralnonlocal} into~\eqref{testinglateral}
we can write  
\begin{equation}\label{testinglateral2}
\begin{split}
\eta^{-\alpha} o_\eta(1) \leq & \ C(r^{-\alpha} + \eta^{-\alpha}) 
+ \I[\mathcal{A}_{\delta}^{ext}](u^\varphi(\cdot, \bar{t}), \bar{x}, D\Psi(\bar{x}, \bar{t}))\\
& - H(\bar{x}, \bar{t}, u(\bar{x}, \bar{t}), D\Psi(\bar{x}, \bar{t})),
\end{split}
\end{equation}
where $D\Psi(\bar{x}, \bar{t}) = \epsilon^{-1} o_\eta(1) + k \eta^{-1} Dd(\bar{x})$. 

Since $u(x_0, t_0) = \varphi(x_0) + \nu$, by the continuity of $\varphi$ and the last fact in~\eqref{proplateral},
for all $\eta$ small enough, using that $u^\varphi=\varphi$ in $\pexterior$ we can write
\begin{equation*}
\begin{split}
\I[\mathcal{A}_{\delta}^{ext}](u^\varphi(\cdot, \bar{t}), \bar{x}, D\Psi(\bar{x}, \bar{t}))
\leq & \ -C \nu \int_{\mathcal{A}_\delta^{ext}} K^\alpha(z)dz \\
& + C(\epsilon^{-1} + \eta^{-1}) \int_{\mathcal{A}_\delta^{ext}} |z| K^\alpha(z)dz. 
\end{split}
\end{equation*}
where we supress the last integral term when $\alpha < 1$. Using the definition of $K^\alpha$, and recalling the choice of $\epsilon$ above,
we conclude from the above inequality that
\begin{equation}\label{nonlocaljump}
\begin{split}
\I[\mathcal{A}_{\delta}^{ext}](u^\varphi(\cdot, \bar{t}), \bar{x}, D\Psi(\bar{x}, \bar{t}))
\leq -C \nu \int_{\mathcal{A}_\delta^{ext}} K^\alpha(z)dz + C\eta^{-1} \tilde{\varrho}_\alpha(\delta),
\end{split}
\end{equation}
where $\tilde{\varrho}(\delta) = \delta^{1 - \alpha}$ if $\alpha > 1$, $\tilde{\varrho}(\delta) = |\ln(\delta)| + 1$ when $\alpha = 1$
and $\tilde{\varrho}(\delta) = 0$ if $\delta < 1$.

At this point we split the analysis. When we consider case $(i.1)$, we just have condition $K$ is nonnegative and bounded, 
and therefore we only can insure that
\begin{equation*}
-C \nu \int_{\mathcal{A}_\delta^{ext}} K^\alpha(z)dz \leq 0.
\end{equation*}

Using this into~\eqref{nonlocaljump} we get 
\begin{equation*}
\I[\mathcal{A}_{\delta}^{ext}](u^\varphi(\cdot, \bar{t}), \bar{x}, D\Psi(\bar{x}, \bar{t})) \leq \tilde{\varrho}(\delta),
\end{equation*}
and replacing this into~\eqref{testinglateral2}, we choose $\delta = \eta$. Applying the definition of $\tilde{\varrho}$ 
and using the condition over $H$ in $(i.1)$, we arrive at
\begin{equation*}
\eta^{-\alpha} o_\eta(1) \leq C(r^{-\alpha} + \eta^{-\alpha})- C_0 k^m \eta^{-m} + \tilde{C},
\end{equation*}
where $\tilde{C}$ depends only on $||u||_\infty$ and the data. We fix $r > 0$ and since $k > 0$ and $m > \alpha$, by choosing $\eta$ small enough,  we reach the contradiction. 

For the case $(i.2)$, recalling that $\alpha < 1$ and the strong ellipticity assumption \eqref{uellipticity},
we have from~\eqref{nonlocaljump} that
$$
\I[\mathcal{A}_{\delta}^{ext}](u^\varphi(\cdot, \bar{t}), \bar{x}, D\Psi(\bar{x}, \bar{t})) \leq -C \nu \delta^{-\alpha}
$$ 
with $C > 0$ independent of $\eta$ and $\delta$. We replace this estimate 
into~\eqref{testinglateral2} to conclude this time that
\begin{equation*}
\eta^{-\alpha}o_\eta(1) \leq C(r^{-\alpha} + \eta^{-\alpha}) - C \nu \delta^{-\alpha} 
- H(\bar{x}, \bar{t}, u(\bar{x}, \bar{t}), k \eta^{-1} Dd(\bar{x}) (1 + o_\eta(1))). 
\end{equation*}

At this point we choose $d(\bar{x}) < \delta < \eta o_\eta(1)$ and applying the condition over the Hamiltonian for this case together 
with~\eqref{proplateral}, we arrive at
\begin{equation*}
\eta^{-\alpha}o_\eta(1) \leq C(r^{-\alpha} + \eta^{-\alpha}) - C \nu \delta^{-\alpha} + c_0 k \eta^{-1} d(\bar{x}) + \tilde{C} ,
\end{equation*}
where $\tilde{C}$ depends only on $||u||_\infty$ and the data. Fixing $r > 0$ and recalling that $\eta^{-1} d(\bar{x}) = o_\eta(1)$, 
we reach the contradiction by choosing $\eta$ small enough. This concludes the proof.
\qed

As a corollary of this lemma we have the following
\begin{prop}\label{lateralboundarycoercive}
Let $\varphi \in C_b(\bar{Q}^{ext})$, $\mathcal{I}$ as in~\eqref{operator} and $H$ with coercive form.
Let $u$ be a bounded viscosity subsolution for the problem~\eqref{cauchy} and let $(x_0, t_0) \in \pboundary$. Then, 
$u(x_0, t_0) \leq \varphi(x_0, t_0)$. In particular, $\tilde{u}$ defined in~\eqref{deftildeuvm} satisfies
$\tilde{u}(x_0, t_0) \leq \varphi(x_0, t_0)$.
\end{prop}

By Remark~\ref{rmksupercoercive}, this result holds since it fits into the case $(i.1)$ in Lemma~\ref{lateralboundary}. 
Concerning the Bellman structure of the problem, we have
\begin{prop}\label{lateralboundaryBellman}
Let $\varphi \in C_b(\bar{Q}^{ext})$, $\alpha < 1$, $\I$ as in~\eqref{operatoralphasmall} satisfying~\eqref{uellipticity}, and
$H$ with Bellman form. Let 
Let $u, v$ be bounded viscosity sub and supersolution for~\eqref{cauchy}, respectively, and
$\tilde{u}$, $\tilde{v}$ as in~\eqref{deftildeuv}. Then
\begin{equation*}
\begin{split}
\tilde{u} \leq u \leq \varphi \leq v \leq \tilde{v} \quad & \mbox{on} \ \Gamma_{out}, \\
\tilde{u} \leq u \leq \varphi \quad & \mbox{on} \ \Gamma.
\end{split}
\end{equation*}
\end{prop}

This result holds since it fits into the cases $(i.2)$ and $(ii)$ in Lemma~\ref{lateralboundary}.


\section{Regularity Issues for Coercive and Bellman Problems.}
\label{regularitysection}

\subsection{Regularity for Coercive Problem.}
We consider the stationary equation 
associated to the coercive version of~\eqref{cauchy}
\begin{eqnarray}\label{stationarym}
\left \{ \begin{array}{rcll} - \I[u] + H_0(x, u, Du) &=& A \quad & in \ \Omega \\ 
u &=& \varphi \quad & in \ \Omega^c, \end{array} \right . 
\end{eqnarray}
where $A > 0$, $\varphi \in C_b(\Omega^c)$, $\I$ is a nonlocal operator of order $\alpha$ with the form~\eqref{operator} or~\eqref{operatoralphasmall}
and $H_0$ defined in (A0) has a coercive form (sub or superlinear). 

As it can be seen in~\cite{BKLT}, the superfractional assumption (A1) makes the gradient term the leading one in equation~\eqref{stationarym}, 
and therefore regularity results can be obtained in an analogous way as in the case of first and second-order equations 
with  coercive Hamiltonians in $Du$ (see~\cite{Barles1},~\cite{Barles-book},~\cite{Capuzzo-Dolcetta-Leoni-Porretta} and references therein).
This regularity result is presented here through the following
\begin{prop}\label{holdersubsol}$\mathrm{(}$\cite{BKLT}$\mathrm{)}$
Let $u$ be a bounded usc viscosity subsolution in $\Omega$ to Equation~\eqref{stationarym}.
Then, there exists a constant $C$ such that, for all $x, y \in \Omega$
\begin{eqnarray*}
|u(x) - u(y)| \leq C|x - y|^{\frac{m - \alpha}{m}} 
\end{eqnarray*}
where $C$ depends on the data, $A$ and $||u^\varphi||_\infty$. In particular, $u$ can be extended up to $\bar{\Omega}$ as a 
H\"older continuous function with H\"older exponent $(m - \alpha)/m$.
\end{prop}


Using this result we can obtain a regularity result for parabolic equations which is sufficient to get the comparison principle.
To do so, we need 
to introduce some notations: for $E \subseteq \R^n$ closed and $g : E \times [0,T] \to \R$ a bounded usc function, 
we define the time sup-convolution of $g$ with parameter $\gamma > 0$ as the function $g^\gamma$ given by 
\begin{equation}\label{ugamma}
g^\gamma(x,t) := \sup \limits_{s \in [0,T]} \{ g(x,s) - \gamma^{-1} (s - t)^2 \}, \quad \mbox{for} \ x \in E, t \in [0,T].
\end{equation}

It is well-known that, for each $\gamma > 0$ and $x \in E$, $t \mapsto g^\gamma(x, t)$ is Lipschitz continuous in $[0,T]$, with 
Lipschitz constant $C_\gamma := 4T \gamma^{-1}$. In addition, if $g \in C(E \times [0,T])$, $g^\gamma \to g$ locally uniformly 
in $E \times [0,T]$ as $\gamma \to 0$.

\begin{lema}\label{lemaugamma}
Let $\varphi \in C_b(\bar{Q}_T^{ext})$, $\I$ as in~\eqref{operator} or~\eqref{operatoralphasmall}
and $H$ with coercive form.  Let $u$ be a bounded viscosity subsolution to problem~\eqref{cauchyT}. 
Then, there exists a constant $a_\gamma > 0$, $a_\gamma \to 0$ as $\gamma \to 0$, such that $u^\gamma$ is a viscosity subsolution 
in $\Omega \times [a_\gamma, T]$ of the problem
\begin{equation*}
\left \{ \begin{array}{rll} \partial_t u^\gamma - \mathcal{I}(u^\gamma) + H(x,t, u^\gamma, Du^\gamma) &=  o_\gamma(1)
\quad & \mbox{in} \ \Omega \times [a_\gamma, T] \\
u^\gamma& = \varphi^\gamma \quad & \mbox{in} \ \Omega^c \times [a_\gamma, T], \end{array} \right .
\end{equation*}
where $o_\gamma(1)$ depends only on the time modulus of continuity of the function $f$ given in (A0).
\end{lema}

\noindent
{\bf \textit{Proof:}} By the upper semicontinuity of $u$, for each $(x,t) \in \bar{Q}_T$ there exists $t_\gamma \in [0,T]$ depending on $x$ and $\gamma$ 
such that
\begin{equation*}
u^\gamma(x,t) = u(x, t_\gamma) - \gamma^{-1} (t - t_\gamma)^2.
\end{equation*}

Since $u$ is bounded, we also have that $|t_\gamma - t| \leq (2 ||u||_{L^\infty(\bar{Q}_T)} \gamma)^{1/2}$ and then we initially 
set $a_\gamma$ as twice this last constant.

We start noting that by applying Proposition~\ref{lateralboundarycoercive}, for each $(x,t) \in \pboundaryT$ we can write
\begin{equation*}
u^\gamma(x,t) \leq \varphi(x, t_\gamma) - \gamma^{-1} (t - t_\gamma)^2 \leq \varphi^\gamma(x, t), 
\end{equation*}
and therefore, the (lateral) boundary condition holds in the classical sense.

Now we address the viscosity inequality in $Q_T$.
Let $(\bar{x},\bar{t}) \in Q_T$ and $\phi$ a smooth test-funtion such that $(\bar{x},\bar{t})$ is a maximum for $u^\gamma - \phi$ in 
$B_{\delta_1}(\bar{x}) \times (\bar{t} - \delta_2, \bar{t} + \delta_2)$ for some $\delta_1, \delta_2 > 0$.
Without loss of generality we can assume $\delta_1 < d(\bar{x})$.

Denote as $\bar{t}_\gamma$ the time attaining the supremum in the definition of $u^\gamma(\bar{x} , \bar{t})$ and
$\tilde{\phi}(x,s) = \phi(x, s + \bar{t} - \bar{t}_\gamma)$. Using the definition of $u^\gamma$ and performing a translation argument in time, we conclude that
\begin{equation*}
u(\bar{x}, \bar{t}_\gamma) - \tilde{\phi}(\bar{x}, \bar{t}_\gamma) \geq u(x,s) - \tilde{\phi}(x,s), 
\quad \mbox{for all} \ (x,s) \in B_{\delta_1}(\bar{x}) \times (\bar{t}_\gamma - \delta_2, \bar{t}_\gamma + \delta_2),
\end{equation*}
which is a testing for $u$ at $(\bar{x}, \bar{t}_\gamma)$ with test-function $\tilde{\phi}$. Applying the viscosity inequality for $u$, we 
can write 
\begin{equation}\label{ineqlemmaugamma1}
E_\delta (u^\varphi, \tilde{\phi}, \bar{x}, \bar{t}_\gamma) \leq 0. 
\end{equation}

Now, using the definition of sup-convolution we have
\begin{equation*}
\begin{split}
u(\bar{x} + z, \bar{t}_\gamma) - \gamma^{-1} (\bar{t}_\gamma - \bar{t})^2 \leq u^\gamma(\bar{x} + z, \bar{t}), & \quad z \in \Omega - \bar{x}, \\ 
\varphi(\bar{x} + z, \bar{t}_\gamma) - \gamma^{-1} (\bar{t}_\gamma - \bar{t})^2 \leq \varphi^\gamma(\bar{x} + z, \bar{t})
, & \quad z \in \Omega^c - \bar{x},
\end{split}
\end{equation*}
meanwhile using that $u^\gamma(\bar{x}, \bar{t}) = u(\bar{x}, \bar{t}_\gamma) - \gamma^{-1}(\bar{t}_\gamma - \bar{t})^2$ we conclude
\begin{equation*}
\I[B_{\delta_1}^c](u(\cdot, \bar{t}_\gamma), \bar{x}, D\tilde{\phi}(\bar{x}, \bar{t}_\gamma)) 
\leq \I[B_{\delta_1}^c]((u^\gamma)^{\varphi^\gamma}(\cdot, \bar{t}), \bar{x}, D\tilde{\phi}(\bar{x}, \bar{t}_\gamma)).
\end{equation*}

Finally, by definition of $\tilde{\phi}$ we have
\begin{equation*}
\partial_t \tilde{\phi}(\bar{x}, \bar{t}_\gamma) = \partial_t \phi(\bar{x}, \bar{t}) \quad \mbox{and} \quad 
D\tilde{\phi}(\bar{x}, \bar{t}_\gamma) = D\phi(\bar{x}, \bar{t}).
\end{equation*}

Using these facts into~\eqref{ineqlemmaugamma1} and using the uniform continuity of $f$, we arrive to the desired viscosity inequality for $u^\gamma$.
\qed

Joining Lemmas~\ref{holdersubsol} and~\ref{lemaugamma} we conclude the following
\begin{lema}\label{regularityugamma}
Let $\varphi \in C_b(\bar{Q}_T^{ext})$, $\I$ as in~\eqref{operator} or~\eqref{operatoralphasmall}, and $H$ with coercive form. 
Then, for all $u$ bounded viscosity subsolution 
to problem~\eqref{cauchyT}, there exists $\gamma_0 > 0$ such that, for all $\gamma \leq \gamma_0$,
$u^\gamma \in C^{1 - \alpha/m, 1}(\Omega \times [a_\gamma, T])$, where $u^\gamma$ is defined in~\eqref{ugamma} and 
$a_\gamma$ is the constant given in Lemma~\ref{lemaugamma}.

Moreover, under the above assumptions, 
$\tilde{u}^\gamma \in C^{1 - \alpha/m, 1}(\bar{\Omega} \times [a_\gamma, T])$, where $\tilde{u}$ is defined in~\eqref{deftildeuvm}.
\end{lema}

\noindent
{\bf \textit{Proof:}} The regularity in $t$ comes from the definition of the sup-convolution. For the H\"older regularity in $x$
the idea is to prove that for each $t \in [a_\gamma, T]$, $x \mapsto u^\gamma(x, t)$ is a viscosity solution to a problem 
like~\eqref{stationarym}. Let $x_0 \in \Omega$, $t_0 \in (a_\gamma, T)$ and $\phi$ a test-function for $u^\gamma(t, \cdot)$ at $x_0$.
For $\epsilon > 0$ small, we incorporate the time variable in the following way
\begin{equation*}
(x,s) \mapsto \Phi(x,s) := u^\gamma(x,s) - \phi(x) - \epsilon^{-1}(s - t_0)^2.
\end{equation*}

The function $\Phi$ being bounded and upper semicontinuous in $\bar{Q}_T$, has a maximum point $(\bar{x},\bar{s}) \in \bar{Q}_T$.
Since $\Phi(\bar{x}, \bar{s}) \geq \Phi(x_0, t_0)$, we have $(\bar{s} - t_0)^2 \leq 2 ||u||_\infty\epsilon$, concluding that 
$\bar{s} \to t_0$ as $\epsilon \to 0$. Then, using the upper semicontinuity of $u^\gamma$, we get $\bar{x} \to x_0$ as $\epsilon \to 0$ too.

Using Lemma~\ref{lemaugamma}, we conclude that
\begin{equation*}
\begin{split}
2\epsilon^{-1}(\bar{s} - t_0) - \I[B_\delta](\phi, \bar{x}) - \I[B_\delta^c](u^\gamma(\bar{s}, \cdot), \bar{x}, D\phi(\bar{x})) & \\
+ H(\bar{x}, \bar{s}, u^\gamma(\bar{x}, \bar{s}), D\phi(\bar{x})) & \leq o_\gamma(1),
\end{split}
\end{equation*}
but we remark that $2\epsilon^{-1}(\bar{s} - t_0) \geq C_\gamma$ because of the Lipschitz continuity of $u^\gamma$ (recall that $2\epsilon^{-1}(\bar{s} - t_0)$ is in the time superdifferential of $u^\gamma$ at $(\bar{x}, \bar{s})$).
Letting $\epsilon \to 0$ and controlling the integral terms by the use of Fatou's Lemma, we conclude that $x \mapsto u^\gamma(t, x)$ is a subsolution
to the problem
\begin{equation*}
-\I(u, x) + H_0(x,u,Du) \leq ||f||_\infty + C_\gamma + o_\gamma(1) \quad \mbox{in} \ \Omega
\end{equation*}
for all $t \in [a_\gamma, T]$. Using Proposition~\ref{holdersubsol}, we conclude the result.

Concerning the last part of the lemma, assume $u = \tilde{u}$. Then, to prove that $u^\gamma \in C^{1 - \alpha/m, 1}(\bar{\Omega} \times [a_\gamma, T])$,
it is sufficient to show that $u^\gamma$ is continuous up to the lateral boundary. In fact, for 
$(x_0, t_0) \in \partial \Omega \times [a_\gamma, T - a_\gamma]$, by definition of $u^\gamma$ and since $u = \tilde{u}$, we can write
\begin{equation*}
u^\gamma(x_0, t_0) = u(x_0, s) - \gamma^{-1}(s - t_0)^2 = u(x_k, s_k) - \gamma^{-1}(s - t_0)^2 + o_k(1)
\end{equation*}
for some $s$ depending on $(x_0, t_0)$, $x_k \to x_0$, $x_k \in \Omega$ and $s_k \to s$. Then
\begin{equation*}
\begin{split}
u^\gamma(x_0, t_0) 
& = u(x_k, s_k) - \gamma^{-1}(s_k - t_0)^2 - \gamma^{-1}o_k(1) \\
& \leq u^\gamma(x_k, t_0) - \gamma^{-1}o_k(1),
\end{split}
\end{equation*}
concluding that 
\begin{equation}\label{limsuplim}
u^\gamma(x_0, t_0) \leq \limsup \limits_{\Omega \ni x \to x_0, t \to t_0} u^\gamma(x,t) = 
\lim \limits_{\Omega \ni x \to x_0, t \to t_0} u^\gamma(x,t), 
\end{equation}
where the last equality comes from
$u^\gamma$ is $C^{1 - \alpha/m, 1}(\Omega \times [a_\gamma, T - a_\gamma])$.

Now, taking $\Omega \ni x_k \to x_0$, we clearly have
$
u^\gamma(x_k, t_0) = u(x_k, s_k) - \gamma^{-1} (s_k - t_0)^2,
$
for some $s_k$ depending on $t_0$ and $x_k$. We see that $(s_k)$ is bounded and therefore it converges to some $\bar{s} \in [0,T]$. Dedefining 
$a_\gamma$ smaller, we have $\bar{s} \in [a_\gamma, T - a_\gamma]$. Now, using the usc of $u$ we have
\begin{equation*}
u^\gamma(x_k, t_0)  \leq u(x_0, \bar{s}) - \gamma^{-1}(\bar{s} - t_0)^2 \leq u^\gamma(x_0, t_0), 
\end{equation*}
from which we get the reverse inequality in~\eqref{limsuplim}. This concludes the proof.
\qed

\subsection{Cone Condition for the Bellman Problem.}

The comfortable H\"older continuity property for subsolutions in the coercive case is hardly available in the Bellman case. However, this property
can be replaced by the weaker ``cone condition'' which is sufficient to apply Soner's argument and to get the desired comparison results, 
see~\cite{Barles-DaLio_On},~\cite{Barles-Rouy},~\cite{usersguide}. 
\begin{prop}\label{conecondition}
Let $\varphi \in C_b(\bar{Q}^{ext})$, $\alpha < 1$, $\I$ as in~\eqref{operatoralphasmall} and $H$ with Bellman form.
Let $u$ be a bounded viscosity subsolution to~\eqref{cauchy} and let $\tilde{u}$ as in~\eqref{deftildeuv}.
Then, for each $(x_0, t_0) \in \Gamma \cup \Gamma_{in}$, there exists $C > 0$ and a sequence $(x_k, t_k) \in Q$ such that, as $k \to \infty$
\begin{eqnarray}\label{conesequence}
\left \{ \begin{array}{l} (x_k, t_k) \to (x_0, t_0); \ \tilde{u}(x_k, t_k) \to \tilde{u}(x_0, t_0), \\ 
|x_k - x_0| \geq Cd(x_k), \\
|t_k - t_0| \geq Cd(x_k). \end{array} \right .
\end{eqnarray}
\end{prop}

We provide the proof of the above cone condition for completeness. However, we note that the results of this section 
are the direct extensions to the parabolic framework of the results presented in~\cite{Topp} and therefore we will omit most of the proofs.

To get Proposition~\ref{conecondition}, we need to introduce notation and give an intermediate result. 
For $x \in \bar{\Omega}$, a function $\phi: \bar{\Omega} \to \R$ bounded, in $C^1(\bar{B}_r(x))$ for some $r > 0$, we define the 
\textsl{censored operator} $\mathcal{I}_\Omega(\phi, x)$ as
\begin{equation*}
\mathcal{I}_{\Omega}(\phi, x) = \int_{\Omega - x} [\phi(x+z) - \phi(x)]K^\alpha(z)dz.
\end{equation*}

Associated to this operator, we have the following proposition
\begin{lema}\label{propU}
Let $\varphi \in C_b(\bar{Q}^{ext})$, $\alpha < 1$, $\I$ as in~\eqref{operatoralphasmall} and $H$ with Bellman form. Let
$u$ be a bounded viscosity subsolution to~\eqref{cauchy} and let
$\tilde{u}$ as in~\eqref{deftildeuv}. Let $(x_0, t_0) \in \pboundary$ and $\beta_0 \in \mathcal{B}$ such that 
\begin{equation}\label{bpointinside0}
b_{\beta_0}(x_0, t_0) \cdot Dd(x_0) \geq c_0 
\end{equation}
for some $c_0 > 0$, and consider the function $U: \bar{Q} \to \R$ defined as
\begin{equation*}\label{U}
U(x,t) = \tilde{u}(x,t) + A  d^{1 - \alpha}(x) 
\end{equation*}

Then, there exists $A, a > 0$ such that $U$ is a viscosity subsolution of the equation
\begin{equation*}\label{eqpropU}
\partial_t u - \I_\Omega(u(\cdot, t)) - b_{\beta_0} \cdot Du = 0 
\quad \mbox{in} \ B_a(x_0) \times (t_0 - a, t_0 + a).
\end{equation*}
\end{lema}

We remark that the notion of viscosity subsolution for censored equations is analogous to the one presented in Definition~\ref{defsol}. 

Using this result, we are in position to prove cone condition.

\medskip

\noindent
{\bf \textit{Proof of Proposition~\ref{conecondition}:}} Note that, if either $x_0 \in \Gamma$ or $x_0 \in \Gamma_{in}$, there exists a control
$\beta_0 \in \mathcal{B}$ satisfying~\eqref{bpointinside0} for some $c_0 > 0$. Thus, denoting $b = b_{\beta_0}$
we can take $r > 0$ small enough such that $b(x,t) \cdot Dd(x) > c_0/2$ for all $x \in \bar{\Omega} \cap \bar{B}_r(x_0)$ and $|t - t_0 | < r$.
After rotation in the $x$ variable and a translation in $(x,t)$, we can assume $t_0 = 0$, $x_0 = 0$ and $Dd(x_0) = e_n$ with $e_n = (0,...,0,1)$, 
implying in particular that $b_n(0,0) > 0$. Finally, denote 
$\h_+ = \{(x', x_n) \in \R^n : x_n > 0\}$ and $A = \bar{\h}_+ \cap \bar{\Omega} \cap \bar{B}_r$.

Recalling the function $U$ defined in Lemma~\ref{propU}, we have this function satisfies the equation 
\begin{equation*}
\partial_t U -\I_\Omega(U(\cdot, t)) - b \cdot DU \leq 0 \quad \mbox{on} \ \bar{A} \times (-r, r).
\end{equation*}

By a simple scaling argument, we conclude the function $(y,s) \mapsto U(\gamma y, \gamma s)$ defined in $\gamma^{-1} (A \times (-r,r))$ 
satisfies the equation
\begin{equation}\label{transportgamma}
\partial_t w -\gamma^{1 - \alpha} \mathcal{I}_{\gamma^{-1}\Omega} (w, y) - b_\gamma (y) \cdot Dw(y) \leq 0 
\quad \mbox{on} \ \gamma^{-1}(A \times (-r,r)),
\end{equation}
where $b_\gamma(y,s) = b(\gamma y, \gamma s)$ for each $(y,s) \in \gamma^{-1}(A \times (-r,r))$. Thus, 
the function $\bar{w} : \bar{\h}_+ \times \R \to \R$ defined as
\begin{eqnarray*}
\bar{w}(x,t) = \limsup \limits_{\gamma \to 0, (z,s) \to (x, t)} U (\gamma z, \gamma s)
\end{eqnarray*}
is a viscosity subsolution for the problem
\begin{eqnarray*}\label{transport}
\partial_t w - b_n(0,0) \frac{\partial w}{\partial y_n} - b'(0,0) \cdot D_{y'}w = 0 \quad \mbox{in} \ \bar{\h}_+ \times \R,
\end{eqnarray*}
by classical arguments in half-relaxed limits applied over the equation~\eqref{transportgamma}. It is worth remark that by Lemma~\ref{propU} 
this equation holds up to the boundary and that $b_n(0,0) > 0$.

The maximal solution for the last transport equation with terminal data $\bar{w}(y', 1, \tau)$ (when we cast $y_n$ as the new ``time'' variable) 
is given by the function
$$
W(y', y_n, s) = \bar{w}(y' - b_{n}(0)^{-1} b'(0)(y_n - 1), 1, s + b_{n}(0)^{-1}(y_n - 1)).
$$ 

Since $W$ is maximal, we have $\bar{w}(y, s) \leq W(y,s)$ when $0 \leq y_n \leq 1$. Now, by definition it is clear that $\bar{w}$ 
is upper semicontinuous and then $\bar{w}(0,0) = U(0,0)$, meanwhile by the upper semicontinuity of $u$ at the boundary and the continuity 
of the distance function we have $\bar{w}(y,s) \leq U(0,0)$ for all $y \in \h_+$. Then, recalling $U(0,0)=\tilde{u}(0,0)$, we conclude that
$$
\tilde{u}(0,0) = \bar{w}(0,0) \leq W(0,0) = \bar{w}(b_n(0)^{-1} b'(0), 1, -b_n(0)^{-1}) \leq \tilde{u}(0,0),
$$
this is $\tilde{u}(0,0) = \bar{w}(x_b,t_b)$, with $x_b = (b_n(0)^{-1} b'(0), 1)$ and $t_b = -b_n(0)^{-1}$. By the very definition of 
$\bar{w}$, we have the existence
of sequences $\gamma_k \to 0$, $t_k \to t_b$, $z_k \to x_b$ such that $(x_k, t_k) := (\gamma_k z_k, \gamma_k t_k)$ satisfies 
$(x_k, t_k) \to (0,0)$ and $\tilde{u}(x_k, t_k) \to \tilde{u}(0,0)$. 

Note that by definition of the sequence $(x_k)_k$ we have 
$
x_k = \gamma_k x_b + o(\gamma_k).
$
Using this, we perform a Taylor expansion on $d(x_k)$, obtaining the existence of a point $\bar{x}_k \in \h_+$ with $\bar{x}_k \to 0$ as 
$k \to \infty$ such that
\begin{equation*}
d(x_k) = Dd(\bar{x}_k) \cdot (\gamma_k x_b + o(\gamma_k)).
\end{equation*}

Hence, since $Dd(0) = e_n$ we conclude $d(x_k) = \gamma_k + o(\gamma_k)$. Thus, using the estimates for $x_k$ and $d(x_k)$ we get that
$
d(x_k) \geq (4|x_b|)^{-1} |x_k|,
$
for all $k$ large enough. Recalling that $x_0 = 0$, we conclude that $(x_k)_k$ is the sequence satisfying~\eqref{conesequence}. Finally, for the 
$t$ variable we have $t_k = \gamma_k t_b + o(\gamma_k)$ and then we get $|t_k| \leq (4|t_b|)^{-1} d(x_k)$ for all $k$ large. 
Recalling $x_0 = 0$ we conclude the result.
\qed

\begin{remark}
It is important to note that, considering~\eqref{conesequence} and its proof, the time and space variables are playing the same 
role regarding the cone's condition property. This fact explains why we cannot weaken the time Lipschitz continuity of $H$ given in
assumption (H).
\end{remark}

Following the same ideas given in Proposition~\ref{conecondition}, it is possible to conclude the cone condition for supersolutions in $\Gamma_{in}$.
\begin{prop}\label{superconecondition}
Let $\varphi \in C_b(\bar{Q}^{ext})$, $\alpha < 1$, $\I$ as in~\eqref{operatoralphasmall} and 
$H$ with Bellman form. Let $v$ a bounded viscosity supersolution to~\eqref{cauchy} and let
$\tilde{v}$ as in~\eqref{deftildeuv}.
Then, for each 
$(x_0,t_0) \in \Gamma_{in}$, there exists a sequence $(x_k, t_k)_k$ of points of $Q$ satisfying~\eqref{conesequence} relative to $\tilde{v}$.
\end{prop}

To get the last proposition, a similar result as Lemma~\ref{propU} is needed for supersolutions. This time we cannot get rid of the 
nonlinearity of $H$ because of the Bellman form, but this can be handled because all the drift terms are pointing ``strictly inside'' $\Omega$.
See~\cite{Topp} for details.


\section{Proof of The Comparison Results.}
\label{proofsection}

\subsection{Strong Comparison Principle for the Coercive Case.} We start with the following
\begin{lema}\label{lemalinearizationcoercive}
Let $\varphi \in C_b(\pexteriorT)$, $\I$ as in~\eqref{operator}, and $H$ with coercive form
satisfying (H1)-(H2).
Let $u, v$ be bounded, respective sub and supersolution to the problem
\begin{eqnarray}\label{eqtolinearize}
\left \{ \begin{array}{rll} \partial_t u - \I(u, x) + H(x,t,u,Du) & = 0 \quad & \mbox{in} \ Q_T \\
u & = \varphi \quad & \mbox{in} \ \pexteriorT, \end{array} \right .
\end{eqnarray}
and let $\tilde{u}$ as in~\eqref{deftildeuvm}.

Let $\gamma \in (0,1)$ and $\mu \in (0,1)$ if $H$ is superlinearly coercive, $\mu = 1$ if $H$ is sublinearly coercive. 
Define $\bar{u} = \mu \tilde{u}^\gamma$ where $u^\gamma$ as in~\eqref{ugamma}, and $w = \bar{u} - v$.
Then, $w$ is a viscosity subsolution for the problem
\begin{equation}\label{eqlinearcoercive}
\left \{ \begin{array}{rll} \partial_t w + h_R w - \I(w, x) - \bar{\omega}_R(|Dw|) & = \bar{C}_R(1 - \mu) + o_\gamma(1) \quad & \mbox{in} \ Q_T \\
w & = \bar{\varphi} - \varphi \quad & \mbox{in} \ \pexteriorT, \end{array} \right .
\end{equation}
where $a_\gamma$ is given in Lemma~\ref{lemaugamma}, $o_\gamma(1)$ depends only on the modulus of continuity of $f$, 
$R = ||\bar{u}||_\infty + ||v||_\infty$, $\bar{\omega}_R$ is a modulus of continuity depending on $R$ and the data, $h_R$
arises in (H1), $\bar{C}_R$ depends on $R$ and $||f||_\infty$, and $\bar{\varphi} = \mu \varphi^\gamma$.
\end{lema}

\noindent
{\bf \textit{Proof:}} We omit the superscript $\sim$ for simplicity and
we address the superlinear case; the sublinear case follows the same ideas with easier computations. 

Note that by Lemma~\ref{lemaugamma} and direct arguments of the viscosity theory, 
we have $\bar{u}$ is a viscosity subsolution to the problem
\begin{equation*}
\begin{split}
\partial_t \bar{u} - \I(\bar{u}, x) + \mu H(x,t, \mu^{-1} \bar{u}, \mu^{-1} D\bar{u}) & = o_\gamma(1) \quad  \mbox{in} \ \Omega \times [a_\gamma, T] \\
\bar{u} & = \bar{\varphi}, \quad \mbox{in} \ \Omega^c \times [a_\gamma, T], 
\end{split} 
\end{equation*}
where $o_\gamma(1) \to 0$ as $\gamma \to 0$ uniformly on $\mu \in (0,1)$. Moreover, by Lemma~\ref{regularityugamma}, we see that
$\bar{u} \in C^{1 -\alpha/m, 1}(\bar{\Omega} \times [a_\gamma, T])$.

The aim is prove that $w$ is a subsolution to~\eqref{eqlinearcoercive} in the viscosity sense with generalized boundary condition, and
the most difficult scenario is when we study the subsolution's obstacle requirement at the lateral boundary.

Let $(x_0, t_0) \in \pboundaryT$. If $w(x_0, t_0) \leq (\bar{\varphi} - \varphi)(x_0, t_0)$, then the boundary condition for 
subsolutions is satisfied in the classical sense and we get the result. For this, we 
assume $w(x_0, t_0) > (\bar{\varphi} - \varphi)(x_0, t_0)$ and the rest of the proof is devoted to conclude the subsolution's viscosity 
inequality at $(x_0, t_0)$. In this case, $w^{\bar{\varphi} - \varphi}(x_0, t_0) = w(x_0, t_0)$, and 
by Lemma~\ref{lateralboundarycoercive} we see that
\begin{equation}\label{vleqvarphi}
v(x_0, t_0) < \bar{u}(x_0, t_0) - \bar{\varphi}(x_0, t_0) + \varphi(x_0, t_0) \leq \varphi(x_0, t_0).
\end{equation}

Let $\phi$ smooth such that $w^{\bar{\varphi} - \varphi} - \phi$ has a strict maximum point in $\bar{Q}_T$ at $(x_0, t_0)$. 
Define $\nu_0 = (Dd(x_0), 0)$ and for all $\epsilon > 0$ we consider the function
\begin{equation*}
\phi_\epsilon(x,y,s,t) = \phi(y,t) + |\epsilon^{-1}((x,s) - (y,t)) - \nu_0|^2. 
\end{equation*}

Now we look for maximum points of the function $\Phi: \bar{\Omega} \times \R^n \times [0,T]^2 \to \R$ defined as
\begin{equation*}
\Phi(x,y,s,t) := \bar{u}(x,s) - v(y,t) - \phi_\epsilon(x,y,s,t).
\end{equation*}

Note that by the boundedness and the upper semicontinuity of $\Phi$, there exists a point 
$(\bar{x}, \bar{y}, \bar{s}, \bar{t}) \in \bar{\Omega} \times \R^n \times [0,T]^2$ attaining the maximum of $\Phi$ in this set.
Then, using the inequality 
$$
\Phi(\bar{x}, \bar{y}, \bar{s}, \bar{t}) \geq \Phi(x_0 + \epsilon Dd(x_0), x_0, t_0, t_0),
$$
together with the continuity of $\bar{u}$ given by Lemma~\ref{regularityugamma}, 
classical arguments in viscosity solution's theory allows us to write
\begin{equation}\label{proplinearization}
\begin{split}
& (\bar{x}, \bar{s}), (\bar{y}, \bar{t}) \to (x_0, t_0), \quad |\epsilon^{-1}(\bar{x} - \bar{y}, \bar{s} - \bar{t}) - \nu_0| \to 0, \\
& \mbox{and} \quad \bar{u}(\bar{x}, \bar{s}) \to \bar{u}(x_0, t_0), \ v_{\varphi}(\bar{y}, \bar{t}) \to v(x_0, t_0), 
\end{split}
\end{equation}
as $\epsilon \to 0$. Moreover, if $\epsilon$ is small enough, we have $\bar{y} \in \bar{\Omega}$, since otherwise, by the 
continuity of $\varphi$, we would have 
\begin{equation*}
v_{\varphi}(\bar{y}, \bar{t}) = \varphi(\bar{x}, \bar{t}) \to \varphi(x_0, t_0)
\end{equation*}
as $\epsilon \to 0$, which is a contradiction to~\eqref{vleqvarphi} in view of the last fact in~\eqref{proplinearization}.
Moreover, by the continuity of $\varphi$ we see that
$v_{\varphi}(\bar{y}, \bar{t}) < \varphi(\bar{y}, \bar{t})$ for all $\epsilon$ small and therefore, 
even if $\bar{y} \in \partial \Omega$, we have a viscosity supersolution inequality associated to $v_{\varphi}$ at $(\bar{y}, \bar{t})$.

On the other hand, by the second property 
in~\eqref{proplinearization} we have
\begin{equation}\label{orderofbarx-bary}
\bar{x} = \bar{y} + \epsilon Dd(x_0) + o_\epsilon(\epsilon),
\end{equation}

A simple Taylor expansion on the distance function implies that $d(\bar{x}) \geq d(\bar{y}) + \epsilon(1 - o_\epsilon(1))$ for all 
$\epsilon$ small enough, concluding that $\bar{x} \in \Omega$. We consider $0 < \delta' < \delta$ and we subtract the viscosity inequality for 
$v$ at $(\bar{y}, \bar{t})$ to the viscosity inequality for $\bar{u}$ at $(\bar{x}, \bar{s})$, concluding that
\begin{equation}\label{testinglinearization}
\mathcal{A} - \I^{\delta'}  \leq o_\gamma(1),
\end{equation}
where 
\begin{equation*}
\begin{split}
\I^{\delta'} = & \ \I[B_{\delta'}](\phi_\epsilon(\cdot, \bar{y}, \bar{s}, \bar{t}), \bar{x}) 
- \I[B_{\delta'}](-\phi_\epsilon(\bar{x}, \cdot, \bar{s}, \bar{t}), \bar{y}) \\
& \ + \I[B_{\delta'}^c](\bar{u}(\cdot, \bar{s}),\bar{x}, \bar{p}) - \I[B_{\delta'}^c](v(\cdot, \bar{t}),\bar{y}, \bar{q}), 
\end{split}
\end{equation*}
with 
\begin{equation*}
\begin{split}
\bar{p} = & \ D_x \phi_\epsilon(\bar{x},\bar{y}, \bar{s}, \bar{t}) = \epsilon^{-1} (\epsilon^{-1}((\bar{x}, \bar{s}) - (\bar{y}, \bar{t})) - \nu_0) \\
\bar{p} = & \ -D_y \phi_\epsilon(\bar{x},\bar{y}, \bar{s}, \bar{t}) = \bar{p} - D\phi(\bar{y}, \bar{t}),
\end{split}
\end{equation*}
and 
\begin{equation*}
\mathcal{A} = (\partial_s \phi_\epsilon - \partial_t \phi_\epsilon)(\bar{x},\bar{y}, \bar{s}, \bar{t})
+ \mu H(\bar{x}, \bar{s}, \mu^{-1} \bar{u}(\bar{x}, \bar{s}), \mu^{-1} \bar{p})- H(\bar{y}, \bar{t}, v(\bar{y}, \bar{t}), \bar{q}).
\end{equation*}

Now we estimate each term in~\eqref{testinglinearization}, starting with $\mathcal{A}$. We have
\begin{equation}\label{Alinearization1}
(\partial_s \phi_\epsilon - \partial_t \phi_\epsilon)(\bar{x},\bar{y}, \bar{s}, \bar{t}) = \partial_t \phi(\bar{y}, \bar{t}),
\end{equation}
and then it remains to estimate the difference among the Hamiltonians to complete the bound for $\mathcal{A}$. 
Using (A0) and the first statement in~\eqref{proplinearization}, we readily have
\begin{equation}\label{Alinearization2}
\mu H(\bar{x}, \bar{s}, \mu^{-1} \bar{u}(\bar{x}, \bar{s}), \mu^{-1} \bar{p})- H(\bar{y}, \bar{t}, v(\bar{y}, \bar{t}), \bar{q})
\geq  (\mu - 1)||f||_\infty - o_\epsilon(1) + \mathcal{H}_0,
\end{equation}
where $o_\epsilon(1) \to 0$ as $\epsilon \to 0$ uniformly in the rest of the variables and $\mathcal{H}_0$ is defined as
\begin{equation*}
\mathcal{H}_0 = \mu H_0(\bar{x}, \mu^{-1} \bar{u}(\bar{x}, \bar{s}), \mu^{-1} \bar{p})- H_0(\bar{y}, v(\bar{y}, \bar{t}), \bar{q}). 
\end{equation*}

Now, using (H1),(A1-b) and (A2-b) we have
\begin{equation*}
\begin{split}
\mathcal{H}_0 \geq & \ h_R(\bar{x}) (\bar{u}(\bar{x}, \bar{s}) - v(\bar{y}, \bar{t})) + (1 - \mu) \Big{(}(m -  1)a_0|\bar{p}|^m - C_R \Big{)} \\
& \ - \omega_R(|\bar{x} - \bar{y}|)(1 + |\bar{p}|^m) - \omega_R(|D_y \phi(\bar{y}, \bar{t})|) |\bar{p}|^{m - 1},
\end{split}
\end{equation*}
where $R = ||\bar{u}||_\infty + ||v||_\infty$. Thus, using the first fact in~\eqref{proplinearization}, 
for all $\epsilon$ small in terms on $1 - \mu$ we can write
\begin{equation*}
\begin{split}
\mathcal{H}_0 \geq & \ (1 - \mu) (m -  1)a_0|\bar{p}|^m/2 - \omega_R(|D_y \phi(\bar{y}, \bar{t})|) |\bar{p}|^{m - 1} \\
& \ + h_R(\bar{x}) (\bar{u}(\bar{x}, \bar{s}) - v(\bar{y}, \bar{t})) - C_R(1 - \mu) - o_\epsilon(1) \\
\geq & \ \inf \limits_{p \geq 0} \{ (1 - \mu) (m -  1)a_0 p^m/2 - \omega_R(|D_y \phi(\bar{y}, \bar{t})|) p^{m - 1} \} \\
& \ + h_R(\bar{x}) (\bar{u}(\bar{x}, \bar{s}) - v(\bar{y}, \bar{t})) - C_R(1 - \mu) - o_\epsilon(1).
\end{split}
\end{equation*}

We notice that the infimum in the last expression is attained, from which we conclude that
\begin{equation*}
\begin{split}
\mathcal{H}_0 \geq & \ - c_{m, \mu} \omega_R(|D\phi(\bar{y}, \bar{t})|)^m \\
& \ + h_R(\bar{x}) (\bar{u}(\bar{x}, \bar{s}) - v(\bar{y}, \bar{t})) - C_R(1 - \mu) - o_\epsilon(1),
\end{split}
\end{equation*}
where $c_{m, \mu} = \frac{(2(m - 1))^{m - 1}}{m^m} ((1 - \mu)(m - 1)a_0)^{1 - m}$. Replacing this into~\eqref{Alinearization2} and 
recalling~\eqref{Alinearization1}, we conclude the following estimate for $\mathcal{A}$
\begin{equation}\label{Alinearization}
\begin{split}
\mathcal{A} \geq & \ \partial_t \phi(\bar{y}, \bar{t}) + h_R(\bar{x}) (\bar{u}(\bar{x}, \bar{s}) - v(\bar{y}, \bar{t}))
- c_{m, \mu} \omega_R(|D\phi(\bar{y}, \bar{t})|)^m \\
& \ + (\mu - 1)(||f||_\infty + C_R) - o_\epsilon(1),  
\end{split}
\end{equation}
where $o_\epsilon(1) \to 0$ as $\epsilon \to 0$ if we keep $\mu, R$ fixed.

Now we addres the estimates for $\I^{\delta'}$. We start noting that
\begin{equation}\label{Ilinearization1}
\I[B_{\delta'}](\phi_\epsilon(\cdot, \bar{y}, \bar{s}, \bar{t}), \bar{x}) 
- \I[B_{\delta'}](-\phi_\epsilon(\bar{x}, \cdot, \bar{s}, \bar{t}), \bar{y}) \leq \epsilon^{-2} o_{\delta'}(1),
\end{equation}
where $o_{\delta'}(1)$ is independent of $\epsilon$. To estimate the integral terms outside $B_{\delta'}$, 
we consider the sets
\begin{equation}\label{defD's}
\begin{split}
D_{int} = (\Omega - \bar{x}) \cap (\Omega - \bar{y}), \quad & D_{ext} =  (\Omega - \bar{x})^c \cap (\Omega - \bar{y})^c, \\ 
D_{int}^{\bar{x}} = (\Omega - \bar{x}) \cap (\Omega - \bar{y})^c, \quad & D_{int}^{\bar{y}} = (\Omega - \bar{x})^c \cap (\Omega - \bar{y}),
\end{split}
\end{equation}
and then we can write
\begin{equation*}
\I[B_{\delta'}^c](\bar{u}(\cdot, \bar{s}),\bar{x}, \bar{p}) - \I[B_{\delta'}^c](v(\cdot, \bar{t}),\bar{y}, \bar{q}) \\
= \I^{\delta'}_{int} + \I^{\delta'}_{int, \bar{x}} + \I^{\delta'}_{int, \bar{y}} + \I^{\delta'}_{ext},
\end{equation*}
where
\begin{equation*}
\begin{split}
\I^{\delta'}_{int} = & \int_{D_{int} \setminus B_{\delta'}} [\bar{u}(\bar{x} + z) - v(\bar{y} + z) - (\bar{u}(\bar{x}) - v(\bar{y})) 
- \mathbf{1}_B \langle D\phi(\bar{y}), z\rangle] K^\alpha(z)dz \\
\I^{\delta'}_{int, \bar{x}} = & \int_{D_{int}^{\bar{x}} \setminus B_{\delta'}} [\bar{u}(\bar{x} + z) - \varphi(\bar{y} + z) 
- (\bar{u}(\bar{x}) - v(\bar{y})) - \mathbf{1}_B \langle D\phi(\bar{y}), z\rangle] K^\alpha(z)dz \\
\I^{\delta'}_{int, \bar{y}} = & \int_{D_{int}^{\bar{y}} \setminus B_{\delta'}} [\bar{\varphi}(\bar{x} + z) - v(\bar{y} + z) 
- (\bar{u}(\bar{x}) - v(\bar{y})) - \mathbf{1}_B \langle D\phi(\bar{y}), z\rangle] K^\alpha(z)dz \\
\I^{\delta'}_{ext} = &\int_{D_{ext} \setminus B_{\delta'}} [\bar{\varphi}(\bar{x} + z) - \varphi(\bar{y} + z) - (\bar{u}(\bar{x}) - v(\bar{y})) 
- \mathbf{1}_B \langle D\phi(\bar{y}), z\rangle] K^\alpha(z)dz.
\end{split}
\end{equation*}

We estimate each integral term separately. For $\I^{\delta'}_{int}$, using that $(\bar{x}, \bar{y}, \bar{s}, \bar{t})$ is a maximum point for $\Phi$
in $\bar{\Omega} \times \R^n \times [0, T]^2$, for all $z \in D_{int}$ we see that
\begin{equation*}
\bar{u}(\bar{x} + z) - v(\bar{y} + z) - (\bar{u}(\bar{x}) - v(\bar{y})) \leq \phi(\bar{y} + z) - \phi(\bar{y}),
\end{equation*}
and therefore we can write
\begin{equation*}
\I^{\delta'}_{int} \leq \I[D_{int} \cap B_\delta \setminus B_{\delta'}](\phi, \bar{y}) + \I^\delta_{int}.
\end{equation*}

We can use the same argument for $I^{\delta'}_{int, \bar{x}}$, concluding that
\begin{equation*}
\I^{\delta'}_{int, \bar{x}} \leq \I[D_{int, \bar{x}} \cap B_\delta \setminus B_{\delta'}](\phi, \bar{y}) + \I^\delta_{int, \bar{x}},
\end{equation*}
but in this case we note that keeping $\delta > 0$ fixed, $ \mathbf{1}_{D_{int, \bar{x}} \setminus B_\delta}(z) K^\alpha(z)$
is an integrable kernel, uniformly in $\epsilon$. Since $|D_{int, \bar{x}}| \to 0$ as $\epsilon \to 0$, we conclude that
\begin{equation*}
\I^{\delta'}_{int, \bar{x}} \leq \I[D_{int, \bar{x}} \cap B_\delta \setminus B_{\delta'}](\phi, \bar{y}) + o_\epsilon(1). 
\end{equation*}

For $\I^{\delta'}_{ext}$, we recall that $w(x_0, t_0) - (\bar{\varphi} - \varphi)(x_0, t_0) > 0$. Then,
by the last fact in~\eqref{proplinearization}, the continuity of $\bar{\varphi}, \varphi$ and the boundedness of $D\phi(\bar{y})$, there exists
$0 < r_0 < \delta$ small not depending on $\epsilon, \delta, \delta'$ such that, for all $r < r_0$ and 
for all $\epsilon$ small enough, we have the inequality
\begin{equation*}
\bar{\varphi}(\bar{x} + z) - \varphi(\bar{y} + z) - (\bar{u}(\bar{x}) - v(\bar{y})) - \langle D\phi(\bar{y}), z\rangle \leq 0, \quad 
\mbox{for all} \ z \in D_{ext} \cap B_{r},
\end{equation*}
and therefore, we arrive at
\begin{equation*}
\I^{\delta'}_{ext} \leq \I^{r}_{ext}. 
\end{equation*}

We finish the estimates for the nonlocal term with $I_{int, \bar{y}}^{\delta'}$. We claim that $D_{int}^{\bar{y}}$ is away from the origin 
uniformly in $\epsilon$ and $\delta'$. This fact is less obvious so we postpone 
its proof until the end. Thus, since $\mathbf{1}_{D_{int} \setminus B_{\delta'}}(z) K\alpha(z)$ is an integrable kernel, 
uniformly in $\delta'$ and $\epsilon$, and since $|D_{int}^{\bar{y}}| \to 0$ as $\epsilon \to 0$, we conclude 
\begin{equation*}
I_{int, \bar{y}}^\delta = o_\epsilon(1). 
\end{equation*}

Thus, joining the above inequalities concerning the integral terms outside $B_{\delta'}$ and~\eqref{Ilinearization1}, we conclude that
\begin{equation*}
\I^{\delta'} \leq \I[(D_{int} \cup D_{int, \bar{x}}) \cap B_\delta \setminus B_{\delta'}](\phi, \bar{y}) + \I^\delta_{int}
+ \I^{r}_{ext} + o_\epsilon(1) + \epsilon^{-2}o_{\delta'}(1),
\end{equation*}
and replacing this and~\eqref{Alinearization} into~\eqref{testinglinearization}, we arrive to
\begin{equation*}
\begin{split}
& \ \partial_t \phi(\bar{y}, \bar{t}) + h_R(\bar{x}) (\bar{u}(\bar{x}, \bar{s}) - v(\bar{y}, \bar{t})) 
- c_{m, \mu} \omega_R(|D\phi(\bar{y}, \bar{t})|)^m \\
& \ - \I[(D_{int} \cup D_{int, \bar{x}}) \cap B_\delta \setminus B_{\delta'}](\phi, \bar{y}) - \I^\delta_{int} - \I^{r}_{ext} \\
\leq & \ (1 - \mu)(||f||_\infty + C_R) + o_\gamma(1) + o_\epsilon(1) + \epsilon^{-2}o_{\delta'}(1).
\end{split}
\end{equation*}

At this point, letting $\delta' \to 0$ and then $\epsilon \to 0$, by~\eqref{proplinearization}, the smoothness of $\phi$,
the continuity of $h_R, \omega_R$ and using Dominated Convergence Theorem, we arrive at
\begin{equation*}
\begin{split}
\partial_t \phi(x_0, t_0) + h_R(x_0) w(x_0, t_0) - c_{m, \mu} \omega_R(|D\phi(x_0, t_0)|)^m &  \\
- \I[(\Omega - x_0) \cap B_\delta](\phi, x_0) & \\
- \I[(\Omega - x_0) \setminus B_\delta](w, x_0, D\phi(x_0)) &  \\
- \I[(\Omega - x_0)^c \setminus B_r](w^{\bar{\varphi} - \varphi}, x_0, D\phi(x_0)) & \leq \bar{C}_R(1 - \mu) + o_\gamma(1),
\end{split}
\end{equation*}
where $\bar{C}_R = ||f||_\infty + C_R$.
Using that $(x_0, t_0)$ is a maximum point for $w^{\bar{\varphi} - \varphi} - \phi$, we can write
\begin{equation*}
\begin{split}
\partial_t \phi(x_0, t_0) + h_R(x_0) w(x_0, t_0) - c_{m, \mu} \omega_R(|D\phi(x_0, t_0)|)^m &  \\
- \I[(\Omega - x_0) \cap B_\delta](\phi, x_0) & \\
- \I[B_\delta^c](w^{\bar{\varphi} - \varphi}, x_0, D\phi(x_0)) &  \\
- \I[(\Omega - x_0)^c \cap B_\delta \setminus B_r](\phi, x_0) & \leq \bar{C}_R(1 - \mu) + o_\gamma(1),
\end{split}
\end{equation*}
and from this, by the smoothness of $\phi$ we can let $r \to 0$, concluding that
\begin{equation*}
\begin{split}
\partial_t \phi(x_0, t_0) + h_R(x_0) w(x_0, t_0) - c_{m, \mu} \omega_R(|D\phi(x_0, t_0)|)^m &  \\
- \I[B_\delta](\phi, x_0) - \I[B_\delta^c](w^{\bar{\varphi} - \varphi}, x_0, D\phi(x_0)) & \leq \bar{C}_R(1 - \mu) + o_\gamma(1),
\end{split}
\end{equation*}
from which we conclude the result.

Now we address the claim leading to the estimate of $\I_{int, \bar{y}}^{\delta'}$. Assume that there exists a sequence $\epsilon_k \to 0$ and 
$z_k \in D_{int}^{\bar{y}}$ such that $z_k \to 0$. 
By definition, there exists
$a_k \in \Omega$ and $b_k \in \Omega^c$ such that $z_k = a_k - \bar{y} = b_k - \bar{x}$ and by the first property in~\eqref{proplinearization} we
have $a_k, b_k \to x_0$. Now, applying~\eqref{orderofbarx-bary} we conclude $b_k = a_k + \epsilon_k(Dd(x_0) + o_{\epsilon_k}(1))$. 
Taking $k$ large we conclude $b_k \in \Omega$, which is a contradiction.
\qed

With the above lemma, we are in position to prove the comparison principle.

\medskip

\noindent
{\bf \textit{Proof of Theorem~\ref{teocoercive}:}} We argue over the redefined function
given by~\eqref{deftildeuvm}, but we omit the superscript $\sim$
for simplicity. We start assuming by contradiction that 
\begin{equation*}
2M := \sup \limits_{\bar{Q}_T} \{ u - v \} > 0.
\end{equation*}

Then, taking $\eta > 0$ small in terms of $M$, we have
\begin{equation}\label{defM}
\sup \limits_{(x,t) \in \bar{Q}_T} \{ u(x,t) - v(x,t) - \eta t\} =: M > 0.
\end{equation}

By the upper semicontinuity of $u - v$ in $\bar{Q}_T$, this supremum is attained at some point $(x_0, t_0) \in \bar{Q}_T$.
By Lemma~\ref{lateralboundarycoercive}, taking $\eta$ smaller if it is necessary, for each $(x_0, t_0)$ attaining $M$ we have $t_0 > 0$.

For the superlinear coercive case, we consider $\eta, \gamma, \mu > 0$, denote $\bar{u} = \mu u^\gamma$ 
and note that $\bar{u} - v - \eta t \to u^\gamma - v$ as $\eta \to 0^+$, $\mu \to 1^-$ uniformly in $\bar{Q}_T$. Since 
$u^\gamma \geq u$ in $\bar{Q}_T$, for all $\eta$ close to 0 and $\mu < 1$ close to 1, we have
\begin{equation}\label{defM2}
\sup \limits_{(x,t) \in \bar{Q}_T} \{ \bar{u}(x,t) - v(x,t) - \eta t\} \geq M/2.
\end{equation}

This supremum is attained at some point $(\tilde{x}, \tilde{t}) \in \bar{Q}_T$. 
Using that $u \leq u^\gamma$, by 
the upper semicontinuity of $u$ and the lower semicontinuity of $v$, we have 
\begin{equation*}
\begin{split}
M \leq & \ \liminf \limits_{\gamma \to 0, \mu \to 1} \{\bar{u}(x_0, t_0) - v(x_0, t_0) - \eta t_0\} \\
\leq & \ \liminf \limits_{\gamma \to 0, \mu \to 1} \{ \sup \limits_{\bar{Q}_T} \{ \bar{u} - v - \eta t\}  \} \\
= & \ \liminf \limits_{\gamma \to 0, \mu \to 1} \{ \bar{u}(\tilde{x}, \tilde{t}) - v(\tilde{x}, \tilde{t}) - \eta \tilde{t}  \} \\
\leq & \ \limsup \limits_{\gamma \to 0, \mu \to 1} \{ \bar{u}(\tilde{x}, \tilde{t}) - v(\tilde{x}, \tilde{t}) - \eta \tilde{t}  \} \leq M,
\end{split}
\end{equation*}
and therefore we have $w(\tilde{x}, \tilde{t}) \to w(x_0, t_0)$ as $\eta, \gamma \to 0$ and $\mu \to 1$, for some $(x_0, t_0)$
attaining $M$ in~\eqref{defM}. In particular, for all $\gamma$ small enough, $\tilde{t} > a_\gamma$, with $a_\gamma$ given in Lemma~\ref{lemaugamma}.

The idea is to use the function $(x,t) \mapsto \eta t$ as test function for $w = \bar{u} - v$
at $(\tilde{x}, \tilde{t})$ and the corresponding viscosity inequality given by Lemma~\ref{lemalinearizationcoercive}.
We can use it at once if $\tilde{x} \in \Omega$ for all $\mu, \gamma$. On the contrary, in the case $\tilde{x} \in \partial \Omega$ we
note that $M/2 \leq w(\tilde{x}, \tilde{t}) = \bar{u}(\tilde{x}, \tilde{t}) - v(\tilde{x}, \tilde{t})$, and by continuity of $\varphi$, we have
$\bar{\varphi} \to \varphi$ locally uniformly in $\Omega^c \times (0,T)$ as $\mu \to 1$ and $\gamma \to 0$. Thus, 
we can take $\mu$ close to 1 and $\gamma$ close to 0 in order to have
\begin{equation*}
w(\tilde{x}, \tilde{t}) > (\bar{\varphi} - \varphi)(\tilde{x}, \tilde{t}),
\end{equation*}
which says that we can test the equation at $(\tilde{x}, \tilde{t})$ even if this point is on the lateral boundary.
Note that this last inequality implies additionally that $w(\tilde{x}, \tilde{t}) = w^{\bar{\varphi} - \varphi}(\tilde{x}, \tilde{t})$.

Thus, for each $\delta > 0$ we can write
\begin{equation*}
\eta + h_R(\tilde{x}) w(\tilde{x}, \tilde{t}) 
- \I[B_\delta^c](w^{\bar{\varphi} - \varphi}(\cdot, \tilde{t}), \tilde{x}, 0) \leq \bar{C}_R (1 - \mu) + o_\gamma(1),
\end{equation*}
where $R = ||\bar{u}||_\infty + ||v||_\infty$. Using that $(\tilde{x}, \tilde{t})$ attains the supremum in~\eqref{defM2} we have
\begin{equation*}
\eta + h_R(\tilde{x}) w(\tilde{x}, \tilde{t}) 
- \I[B_\delta^c \cap (\Omega^c - \tilde{x})](w^{\bar{\varphi} - \varphi}(\cdot, \tilde{t}), \tilde{x}, 0) \leq \bar{C}_R (1 - \mu) + o_\gamma(1), 
\end{equation*}
and from this we see that
\begin{equation*}
\begin{split}
\eta + \Big{(} h_R(\tilde{x}) + \int \limits_{(\Omega^c - \tilde{x}) \setminus B_\delta}K^\alpha(z)dz \Big {)}w(\tilde{x}, \tilde{t}) & \\
- \int \limits_{(\Omega^c - \tilde{x}) \setminus B_\delta}(\bar{\varphi}(\tilde{x} + z) - \varphi(\tilde{x} + z))K^\alpha(z)dz 
\leq & \ \bar{C}_R (1 - \mu) + o_\gamma(1),  
\end{split}
\end{equation*}

But using that $\bar{\varphi} \to \varphi$ locally uniform in $\Omega^c \times (0,T)$ as $\mu \to 1$ and $\eta \to 0$, 
using Dominated Convergence Theorem, the continuity of $h_R$ and that 
$w(\tilde{x}, \tilde{t}) \to M$, taking $\eta, \gamma \to 0$ and  $\mu \to 1$ we arrive at
\begin{equation*}
\eta + \Big{(} h_R(x_0) + \int \limits_{(\Omega^c - x_0) \setminus B_\delta}K^\alpha(z)dz \Big {)}w(x_0, t_0) \leq 0,
\end{equation*}
where $(x_0, t_0)$ is a point attaining the supremum in~\eqref{defM}.
Finally, by (H1) we can take $\delta > 0$ small in order to have
\begin{equation*}
\eta/2 \leq 0, 
\end{equation*}
which is a contradiction.
\qed

\subsection{Strong Comparison Principle for the Bellman Case.}
The analogous to Lemma~\ref{lemalinearizationcoercive} for the Bellman case reads as follows
\begin{lema}\label{lemalinearizationBellman}
Let $\varphi \in C_b(\pexteriorT)$, $\alpha < 1$, $\I$ as in~\eqref{operatoralphasmall}, with $K$ satisfying (UE) and
$H$ with Bellman form, satisfying (H1)-(H2).
Let $u, v$ be bounded, respective viscosity sub and supersolution to~\eqref{eqtolinearize}, and consider $\tilde{u}, \tilde{v}$
as in~\eqref{deftildeuv}.
Then, $w := \tilde{u} - \tilde{v}$ is a viscosity subsolution for the problem
\begin{equation*}
\begin{split}
\partial_t w + h_R(x) w - \I(w, x) - \beta|Dw| & = 0 \quad  \mbox{in} \ Q_T \\
w & = 0 \quad \mbox{in} \ \pexteriorT 
\end{split} 
\end{equation*}
where 
$
\beta = \sup_{\beta \in \mathcal{B}} |b_\beta(x_0, t_0)|,
$
$R = ||\bar{u}||_\infty + ||v||_\infty$, $\bar{\omega}$ is a modulus of continuity depending on $b$ and $h_R$
arises in (H1).
\end{lema}

We require the following result which states the viscosity inequality holds on $\Gamma_{in}$ for the redefined functions $\tilde{u}, \tilde{v}$.
\begin{lema}\label{lemaineqGammain}
Assume the conditions of Lemma~\ref{lemalinearizationBellman} hold. Let $(x_0, t_0) \in \Gamma_{in}$ and 
assume $\tilde{u}(x_0, t_0) > \varphi(x_0, t_0)$. Then, for each $\phi$ smooth such that $(x_0, t_0)$ is a maximum point 
for $\tilde{u}^\varphi - \phi$ in $B_\delta(x_0) \times (t_0 - \delta, t_0 + \delta)$ for some $\delta > 0$, then 
$E_\delta(\tilde{u}^\varphi, \phi, x_0, t_0) \leq 0$. The analogous result holds for $\tilde{v}$.
\end{lema}

\noindent
{\bf \textit{Proof:}} Let $(x_k, t_k) \to (x_0, t_0)$ such that $\tilde{u}(x_k, t_k) \to \tilde{u}(x_0, t_0)$, with $x_k \in \Omega$.
Define $\epsilon_k = d(x_k)$ and consider the function
\begin{equation*}
(x,t) \mapsto \tilde{u}^\varphi(x,t) - \phi(x,t) + \epsilon_k \ln(d(x)) \mathbf{1}_\Omega(x). 
\end{equation*}

For $k$ large enough, we have this function has a maximum point $(\bar{x}_k, \bar{t}_k)$ in $B_\delta(x_0) \times (t_0 - \delta, t_0 + \delta)$,
with $(\bar{x}_k, \bar{t}_k) \to (x_0, t_0)$, $\tilde{u}(\bar{x}_k, \bar{t}_k) \to \tilde{u}(x_0, t_0)$ and $\bar{x}_k \in \Omega$. Using this
and since $u^\varphi = \tilde{u}^\varphi$ up to a set of zero Lebesgue measure, we can write the viscosity inequality for 
$u$ at $(\bar{x}_k, \bar{t}_k)$
\begin{equation*}
E_\delta(\tilde{u}^\varphi, \phi - \epsilon_k \ln(d) \mathbf{1}_\Omega, \bar{x}_k, \bar{t}_k) \leq 0.
\end{equation*}

But using that $(x_0, t_0) \in \Gamma_{in}$, there exists $c_0 > 0$ such that, 
for all $k$ large enough we have $b_\beta(\bar{x}_k, \bar{t}_j) \cdot Dd(\bar{x}_k) \geq c_0$. Thus, we arrive at
\begin{equation*}
\epsilon_k \Big{(} -\I[B_\delta \cap (\Omega - \bar{x}_k)](\ln(d), \bar{x}_k) + c_0 d^{-1}(\bar{x}_k) \Big{)} +
E_\delta(\tilde{u}^\varphi, \phi, \bar{x}_k, \bar{t}_k) \leq 0.
\end{equation*}

Here we mention that there exists a constant $c > 0$ such that 
$$
\I[B_\delta \cap (\Omega - \bar{x}_k)](\ln(d), \bar{x}_k) \leq c d^{-\alpha}(\bar{x}_k),
$$
see~\cite{Topp} for a proof of this result. Thus, for all $k$ large we have
\begin{equation*}
E_\delta(\tilde{u}^\varphi, \phi, \bar{x}_k, \bar{t}_k) \leq 0, 
\end{equation*}
and recalling that $\tilde{u}^\varphi = \tilde{u}$ in a neighborhood of $(x_0, t_0)$, taking $k \to \infty$ together with Dominated Convergence
Theorem to control the integral terms, we get the result.
\qed

\medskip
\noindent
{\bf \textit{Proof of Lemma~\ref{lemalinearizationBellman}:}} We concentrate in the viscosity inequality on the lateral boundary.
By Lemma~\ref{lateralboundaryBellman}, the interesting case is when the test point $(x_0, t_0) \in \Gamma \cup \Gamma_{in}$ is such that
$w(x_0, t_0) > 0$. Note that $w^0(x_0, t_0) = w(x_0, t_0)$ in this case.

Consider $\phi$ a smooth function such that $w^0 -\phi$ has a strict maximum point in $\bar{Q}_T$ at $(x_0, t_0)$.

If $(x_0, t_0) \in \Gamma$, 
Proposition~\ref{lateralboundaryBellman} allows us to conclude $\tilde{u}(x_0, t_0) \leq \varphi(x_0, t_0)$ and Proposition~\ref{conecondition} implies 
the existence of a sequence satisfying~\eqref{conesequence}. In particular, denoting
$
\epsilon_k = \sqrt{|x_k - x_0|^2 + (t_k - t_0)^2} ,
$
up to a subsequences we have $\epsilon_k^{-1}(x_k, t_k) \to \nu_0$ satisfying $\nu_0 \cdot (Dd(x_0), 0) \geq c_0$, for some
$c_0 > 0$. This time, for $k \in \N$ we double variables and use the penalization
\begin{equation*}
\tilde{u}(x,s)  - \tilde{v}(y,t) -  \phi(y, t) - |\epsilon_k^{-1}((x,s) - (y,t)) - \nu_0|^2,
\end{equation*}
and from this point we argue exactly as in Lemma~\ref{lemalinearizationcoercive}, arriving at inequality~\eqref{testinglinearization},
where $\I^{\delta'}$ is managed in the same way as in the coercive case, but $\mathcal{A}$ in this case has the form
\begin{equation*}
\mathcal{A} \geq \partial_t \phi(x_0, t_0) + h_R(x_0) w(x_0, t_0) - \beta |D\phi(x_0, t_0)| - o_k(1),
\end{equation*}
where $o_k(1) \to 0$ as $k \to \infty$.
From this, we proceed exactly as in the proof of Lemma~\ref{lemalinearizationcoercive} to conclude the result.

If $(x_0, t_0) \in \Gamma_{in}$, we consider two sub-cases: if $\tilde{v}(x_0, t_0) < \varphi(x_0, t_0)$, then we argue exactly as
in the case of $\Gamma$ because cone condition also holds for subsolutions on $\Gamma_{in}$. On the other hand, if 
$\varphi(x_0, t_0) \leq \tilde{v}(x_0, t_0)$, we can exchange the roles of $u$ and $v$ in the proof of the case $(x_0, t_0) \in \Gamma$ since 
cone condition holds for supersolution on $\Gamma_{in}$ as it is stated in Proposition~\ref{superconecondition}. We remark that by 
Lemma~\ref{lemaineqGammain} we can use the viscosity inequality on $\Gamma_{in}$ for $\tilde{u}$ and/or $\tilde{v}$ if they do not satisfy the 
boundary condition in the classical sense.
\qed

\medskip
\noindent
{\bf \textit{Proof of Theorem~\ref{teoBellman}:}} We argue by contradiction as in the proof of Theorem~\ref{teocoercive}, where this time
the linearization precedure is played by Lemma~\ref{lemalinearizationBellman}. We omit the details.
\qed


\section{Existence and Large Time Behavior.}
\label{existencesection}

\subsection{Existence and Uniqueness Issues.}
For both coercive and Bellman case, the application of Perron's method on a sequence of finite-time horizon problems 
with the form~\eqref{cauchyT} with $T \to \infty$ and the strong comparison principle allows us to get the existence of a solution 
which is defined for all time.

For reasons that will be made clear in the next theorem , we introduce the following nondegeneracy condition:

\medskip
\noindent
(H2') \textsl{There exists $\mu_0 > 0$ and a continuous function 
$h: \bar{\Omega} \to \R$ satisfying
\begin{equation*}
\inf \limits_{x \in \bar{\Omega}} \Big{\{} h(x) + \int_{x + z \notin \Omega} K^\alpha(z)dz \Big{\}} \geq \mu_0,
\end{equation*}
such that, for all $R > 0$, $h_R$ defined in (H1) satisfies $h_R \geq h$.}

\begin{teo}\label{teoexistence}{\bf (Existence and Uniqueness)}
Let $\alpha \in (0,2)$, $u_0 \in C(\bar{\Omega})$, $\varphi \in C_b(\bar{Q}^{ext})$ satisfying~(H0). Assume~\eqref{cauchy} has a

\medskip
\noindent
$\bullet$ Coercive Form: $\mathcal{I}$ as in~\eqref{operator} (as in~\eqref{operatoralphasmall} if $\alpha < 1$), and $H$ has coercive form. 

\medskip
\noindent
$\bullet$ Bellman Form: $\alpha < 1$, $\I$ as in~\eqref{operatoralphasmall} satisfying (UE), and $H$ has Bellman form.

\medskip

In both cases, we further assume that $H$ satisfies (H1)-(H2).
Then, there exists a unique viscosity solution $u \in C(\bar{Q}) \cap L^\infty(\bar{Q}_T)$ for all $T > 0$, to problem~\eqref{cauchy}.

Moreover, if (H2') holds, then the unique solution $u \in C(\bar{Q}) \cap L^\infty(\bar{Q}_T)$ for all $T > 0$, 
to problem~\eqref{cauchy}, is uniformly bounded in $\bar{Q}$.
\end{teo}

Theorem~\ref{teoexistence} for the finite time horizon problem~\eqref{cauchyT} follows from the application of Perron's method over 
an extended problem over $\R^n \times [0,T]$. For this auxiliary problem, the role of the global sub and supersolution present in Perron's method 
is played by functions with the form $(x,t) \mapsto C_1t + C_2$, for suitable constants $C_1, C_2$ depending on the data and $T$. On the other hand, 
under the assumption (H2') these global sub and supersolution can be taken as constant functions depending on the data, but not on $T$, concluding
the uniform boundedness. See~\cite{Barles-Chasseigne-Imbert},~\cite{Topp} for details.

Assumption (H2') also allows us to get the strong 
comparison principle and therefore the existence and uniqueness for the associated stationary problem.
\begin{teo}\label{teoexistencestationary}
Let $\bar{\varphi} \in C_b(\Omega^c)$, $\bar{H} \in C(\bar{\Omega} \times \R \times \R^n)$ and consider
\begin{eqnarray}\label{inftyeq}
\left \{ \begin{array}{rcll} - \mathcal{I}(u) + \bar{H}(x, u, Du) &=& 0 \quad & in \ \Omega, \\ 
u &=& \bar{\varphi} \quad & in \ \Omega^c. \end{array} \right . 
\end{eqnarray}

Assume this problem has coercive or Bellman form in the sense of Theorem~\ref{teoexistence} in the time independent framework, 
with $\bar{H}$ satisfying (H1), (H2) and (H2'). Then, there exists a unique viscosity solution $u \in C(\bar{\Omega})$ for~\eqref{inftyeq}.
\end{teo}

\subsection{Large Time Behavior.}

Once the existence and uniqueness for problem~\eqref{cauchy} is obtained, it arises the natural question of the asymptotic behavior of the solution 
as $t \to +\infty$. For our models, the answer is contained in the following
\begin{teo}\label{teoconvergence}
Let $u_0 \in C(\bar{\Omega})$ and $\varphi \in C_b(\bar{Q}^{ext})$ satisfying (H0). 
Assume~\eqref{cauchy} has coercive or Bellman form in the sense of Theorem~\ref{teoexistence}, with $H$ satisfying (H1),(H2) and (H2').
Assume there exist continuous functions $\bar{H}: \bar{\Omega} \times \R \times \R^n \to \R$ and $\bar{\varphi} : \Omega^c \to \R$ satisfying
\begin{equation}\label{dataconvergence}
\begin{split}
& H(\cdot, t, \cdot, \cdot) \to \bar{H} \quad \mbox{in} \ C (\bar{\Omega} \times \R \times \R^n), \\
& \varphi(\cdot, t) \to \bar{\varphi} \quad \mbox{in} \ C(\Omega^c),
\end{split}
\end{equation}
as $t \to \infty$. Then, the unique viscosity solution $u$ of~\eqref{cauchy} converges uniformly in $\bar{\Omega}$ to $u_\infty$, the 
unique viscosity solution of the problem~\eqref{inftyeq}.
\end{teo}

\noindent
{\bf \textit{Proof:}} The proof of this theorem can be framed in the general context of parabolic equations for which the limit problem 
satisfied the comparison principle. 
For each $(x,t) \in \bar{\Omega} \times [0,+\infty)$, define the functions
\begin{equation*}
\begin{split}
& \bar{u}(x,t) = \limsup \limits_{\epsilon \to 0, z \to x, z \in \Omega} u(z, t/\epsilon), \\
& \underline{u}(x,t) = \liminf \limits_{\epsilon \to 0, z \to x, z \in \Omega} u(z, t/\epsilon),
\end{split}
\end{equation*}
which are well defined by the uniform boundedness of $u$.
The application of the half-relaxed limits method proves that for all $t > 0$, the functions 
$x \mapsto \bar{u}(x, t)$ and $x \mapsto \underline{u}(x, t)$ are respectively viscosity sub and supersolution for problem~\eqref{inftyeq}.
Then, by comparison principle for Dirichlet problems we have $\bar{u} = \underline{u}$ 
in $\bar{Q}$ and consequently
$\bar{u}(t, x) = \underline{u}(t,x) = u_\infty(x)$ for all $(x,t) \in \bar{Q}$ by the uniqueness of problem~\eqref{inftyeq}. This concludes the result. 
\qed

We can provide a rate of convergence in the particular case that $H$ is time independent and $\varphi$ converges uniformly 
to $\bar{\varphi}$ as $t\to \infty$.
\begin{prop}\label{convergence}
Let $u_0 \in C(\bar{\Omega})$, $\varphi \in C_b(\bar{Q}^{ext})$ satisfying~(H0), and $H \in C(\bar{\Omega} \times \R \times \R^n)$
is time independent.
Assume problem~\eqref{cauchy} has coercive or Bellman form in the sense of Theorem~\ref{teoexistence}, with $H$ satisfying
(H1), (H2) and (H2').
Assume there exists $\bar{\varphi} \in C_b(\Omega^c)$ such that 
$$
\varphi(\cdot, t) \to \bar{\varphi}
$$ 
uniformly in $\Omega^c$ as $t \to \infty$. 

Let $u$ be the unique solution to problem~\eqref{cauchy}, and
$u_\infty$ be the unique bounded viscosity solution to~\eqref{inftyeq} asociated to $\bar{H} = H$ and $\bar{\varphi}$. Then,
\begin{equation*}
||u(\cdot,t) - u_\infty||_{L^\infty(\bar{\Omega})} \leq e^{-\mu_0 t} \Big{(} ||u_0 - u_\infty||_{L^\infty(\bar{\Omega})} 
+ \mu_0 \int_{-\infty}^{t} g (s) e^{\mu_0 s } ds \Big{)},
\end{equation*}
where $g$ is defined as
\begin{equation*}
g(t) = \left \{ \begin{array}{ll} \sup \limits_{\tau \geq t} ||\varphi(\cdot, \tau) - \bar{\varphi}||_{L^\infty(\Omega^c)} 
& \quad \mbox{for} \ t \geq 0 \\
\sup \limits_{\tau \geq 0} ||\varphi(\cdot, \tau) - \bar{\varphi}||_{L^\infty(\Omega^c)} & \quad \mbox{for} \ t < 0.
\end{array} \right .
\end{equation*}
\end{prop}

\noindent
{\bf \textit{Proof:}} Note that $g(t) \leq ||\varphi||_\infty + ||\bar{\varphi}||_\infty$ and then, the function
\begin{equation*}
G(t) = \mu_0 \int_{-\infty}^{t} g(s) e^{\mu_0 s}ds 
\end{equation*}
is well defined. Note also that $g$ is decreasing in $t$ and this implies that
\begin{equation}\label{orderG}
e^{-\mu_0 t} G(t) \geq \mu_0 e^{-\mu_0 t}  g(t) \int_{-\infty}^{t} e^{\mu_0 s}ds \geq g(t).
\end{equation}

With this, consider the function
\begin{equation*}
U(x,t) = u_\infty(x)  + e^{-\mu_0 t} \tilde{G}(t).
\end{equation*}
where $\tilde{G}(t) = G(t) + ||u_0 - u_\infty||_{L^\infty(\bar{\Omega})}$.
We claim $U$ is a supersolution for the problem satisfied by $u$. In fact, for all $x \in \bar{\Omega}$ we clearly have
\begin{equation*}
U(x,0) \geq u_\infty(x) + ||u_0 - u_\infty||_{L^\infty(\bar{\Omega})} \geq u_0(x).
\end{equation*}

Let $(x_0,t_0) \in Q$ and let $\phi$ be a smooth function such that 
$(x_0, t_0)$ is a minimum point of $U_\varphi - \phi$ in $B_\delta(x_0) \times (t_0- \delta, t_0 + \delta)$. 
At one hand, from this testing we have
\begin{equation}\label{ratederivative}
\partial_t \phi(x_0, t_0) = -\mu_0 e^{-\mu_0 t_0} \tilde{G}(t_0) + \mu_0 g(t_0).
\end{equation}

On the other hand, we get that $x_0$ is a minimum point for the function
\begin{equation*}
x \mapsto (u_\infty)_{\bar{\varphi}}(x) - \Big{(} -e^{-\mu_0 t_0} \tilde{G}(t_0) + \phi(x,t_0)\Big{)}
\end{equation*}
in $B_\delta(x_0)$. Hence, we use this as a testing for $u_\infty$, which is a supersolution for the problem~\eqref{inftyeq} at 
$x_0$. Using the viscosity inequality for $u_\infty$, the definition of $U$, the equality~\eqref{ratederivative} and 
the assumption (H1), we arrive to
\begin{equation}\label{ineqconvergence}
\begin{split}
& \partial_t \phi(x_0, t_0) - \I[B_\delta](\phi(\cdot, t_0), x_0) - \I[B_\delta^c](U^\varphi(\cdot, t_0), x_0, D\phi(x_0, t_0)) \\
\geq & - H(x_0, U(x_0, t_0), D\phi(x_0, t_0)) + A_0,
\end{split}
\end{equation}
where
\begin{equation*}
\begin{split}
A_0 = & -\mu_0 e^{-\mu_0 t_0} \tilde{G}(t_0) + \mu_0 g(t_0) + h(x_0)e^{-\mu_0 t_0} \tilde{G}(t_0) \\
& + \int_{(\Omega - x_0)^c \setminus B_\delta} [e^{-\mu_0 t_0} \tilde{G}(t_0) - (\varphi(x_0 + z, t_0) - \bar{\varphi}(x_0 + z))]K^\alpha(z)dz. 
\end{split}
\end{equation*}

But clearly we have
\begin{equation*}
A_0 \geq (e^{-\mu_0 t_0} \tilde{G}(t_0) - g(t_0)) \Big{(} \int_{(\Omega- x_0)^c \setminus B_\delta} K^\alpha(z)dz + h(x_0) - \mu_0 \Big{)},
\end{equation*}
and applying (H2') and~\eqref{orderG}, we obtain $A_0 \geq 0$. This concludes the claim when $(x_0, t_0) \in Q$.
For $(x,t) \in \pboundary$ and $U(x,t) < \varphi(x,t)$, by definition we have
\begin{equation*}
u_\infty(x) < \varphi(x, t) - e^{-\mu_0 t} G(t). 
\end{equation*}

Using the inequality~\eqref{orderG} and the definition of $g$, we conclude
\begin{equation*}
u_\infty(x) < \bar{\varphi}(x), 
\end{equation*}
concluding that in this case we can use the corresponding viscosity inequality for $u_\infty$, concluding the claim.

In the same way a subsolution can be constructed, and the result follows by comparison principle.
\qed

\bigskip

\noindent
{\bf Aknowledgements:} G.B. is partially supported by the ANR (Agence Nationale de
la Recherche) through ANR WKBHJ (ANR-12-BS01-0020). 
E.T. was partially supported by CONICYT, Grants Capital Humano Avanzado, Realizaci\'on de Tesis Doctoral and Cotutela en el 
Extranjero.

\end{document}